\documentclass[11pt]{article}
\usepackage[leqno]{amsmath}
\usepackage{amsfonts,amssymb,theorem,amscd}
\usepackage[mathscr]{eucal}
      \makeatletter
\def\@@seccntfont{\bfseries\slshape}
\def\@@secheaderfont{\bfseries\upshape}
\def\@@precnt{{\upshape}}
\def\@@postcnt{{\upshape}}
\def\@@startsection#1#2#3#4#5#6{\if@noskipsec \leavevmode \fi\par
  \@tempskipa #4\relax\@afterindenttrue
  \ifdim \@tempskipa <\z@\@tempskipa -\@tempskipa \@afterindentfalse\fi
  \if@nobreak\everypar{}\else\addpenalty\@secpenalty\addvspace\@tempskipa\fi
  \@dblarg{\@@sect{#1}{#2}{#3}{#4}{#5}{#6}}}
\newskip\@@argfive
\def\@@sect#1#2#3#4#5#6[#7]#8{\ifnum #2>\c@secnumdepth\let\@svsec\@empty
  \else\refstepcounter{#1}\protected@edef\@svsec{\@seccntformat{#1}\relax}\fi\@tempskipa #5\relax
  \ifdim \@tempskipa>\z@\begingroup#6{\@hangfrom{\hskip#3\relax\@svsec}\interlinepenalty 
  \@M #8\@@par}\endgroup\csname #1mark\endcsname{#7}\else
  \def\@svsechd{#6{\hskip #3\relax\@svsec#8}\csname#1mark\endcsname{#7}}\fi\@xsect{#5}}
\def\@@seccntfmt{\renewcommand{\@seccntformat}[1]{%
  {\@@seccntfont\@@precnt\csname the##1\endcsname}\@@postcnt\hspace{1ex}}}
\newcommand{\@@secnostar}[1][]{\def\@@tmpa{}\def\@@tmpb{#1}%
  \ifx\@@tmpa\@@tmpb\@@argfive-1ex\else\@@argfive-\labelsep\fi\@@seccntfmt\@@startsection
  {\@@name}{\@@level}{0mm}{-\baselineskip}{\@@argfive}{\@@secheaderfont}{#1}}
\newcommand{\@@secstar}[2][]{\def\@@tmpa{}\def\@@tmpb{#1}%
  \ifx\@@tmpa\@@tmpb\def\@@tmpa{#2}\else\def\@@tmpa{#1}\fi\@@seccntfmt\@startsection
  {\@@name}{\@@level}{0mm}{-\baselineskip}{-\labelsep}{\@@secheaderfont}[\@@tmpa]{#2}}
\def\@@section{\def\@@name{section}\def\@@level{1}\@ifstar{\@@secstar}{\@@secnostar}}
\def\@@subsection{\def\@@name{subsection}\def\@@level{2}\@ifstar{\@@secstar}{\@@secnostar}}
\def\@@subsubsection{\def\@@name{subsubsection}\def\@@level{3}\@ifstar{\@@secstar}{\@@secnostar}}
\def\@@paragraph{\def\@@name{paragraph}\def\@@level{4}\@ifstar{\@@secstar}{\@@secnostar}}
\def\@@subparagraph{\def\@@name{subparagraph}\def\@@level{5}\@ifstar{\@@secstar}{\@@secnostar}}
\let\@@latexsection=\section\let\section=\@@section
\let\@@latexsubsection=\subsection\let\subsection=\@@subsection
\let\@@latexsubsubsection=\subsubsection\let\subsubsection=\@@subsubsection
\let\@@latexparagraph=\paragraph\let\paragraph=\@@paragraph
\let\@@latexsubparagraph=\subparagraph\let\subparagraph=\@@subparagraph
\def\setseccntfmt{\renewcommand{\@seccntformat}[1]{\S
  \csname the##1\endcsname\hspace{1ex}}}
\renewcommand{\section}{\setseccntfmt\@startsection
  {section}{1}{0mm}{-\baselineskip}{1.0\baselineskip}{\sf\bfseries\Large}}
\newif\if@@lastcharstar
\def\@@xxlastcharstar#1{\gdef\@@prevchar{}\@@lastcharstarfalse\@@yylastcharstar#1\end}
\def\@@yylastcharstar#1{\ifx#1\end\def\@@tmpa{*}\ifx\@@prevchar\@@tmpa\@@lastcharstartrue\fi
  \let\@@next=\relax\else\def\@@prevchar{#1}\let\@@next=\@@yylastcharstar\fi\@@next}
\gdef\@thm#1#2{\@@xxlastcharstar{#1}\if@@lastcharstar\else\refstepcounter{#1}\fi
  \trivlist\@topsep \theorempreskipamount\@topsepadd \theorempostskipamount
  \@ifnextchar [{\@ythm{#1}{#2}}{\@begintheorem{#2}{\csname the#1\endcsname}\ignorespaces}}
\gdef\th@nonum{\normalfont\itshape
  \def\@begintheorem##1##2{\item[\hskip\labelsep\theorem@headerfont ##1]}%
  \def\@opargbegintheorem##1##2##3{%
    \item[\hskip\labelsep \theorem@headerfont ##1\ (##3)]}}
\gdef\th@change{
  \def\@begintheorem##1##2{\item[\hskip\labelsep{\@@seccntfont
  \@@precnt##2\@@postcnt\hskip 1ex}\theorem@headerfont ##1]}%
  \def\@opargbegintheorem##1##2##3{\item[\hskip\labelsep{\@@seccntfont
  \@@precnt##2\@@postcnt\hskip 1ex}\theorem@headerfont ##1\ (##3)]}}
\def\newtheoremset[#1][#2]{{\theoremstyle{change}\newtheorem{-#1}[section]{#2}%
  \newtheorem{#1}[subsection]{#2}\newtheorem{*#1}[subsubsection]{#2}%
  \newtheorem{**#1}[paragraph]{#2}\newtheorem{***#1}[subparagraph]{#2}}
  {\theoremstyle{nonum}\newtheorem{#1*}{#2}}}
\theoremheaderfont{\normalfont\bfseries}
{\theorembodyfont{\normalfont\itshape}
\newtheoremset[thm][Theorem]\newtheoremset[prop][Proposition]
\newtheoremset[lemma][Lemma]\newtheoremset[sublemma][Sublemma]
\newtheoremset[cor][Corollary]\newtheoremset[conj][Conjecture]}
{\theorembodyfont{\normalfont\rmfamily}
\newtheoremset[exam][Example]
\newtheoremset[defn][Definition]\newtheoremset[rem][Remark]
\newtheoremset[nota][Notation]\newtheoremset[hypo][Hypothesis]}

\def\makeop#1{\expandafter\def\csname#1\endcsname{\mathop{\rm #1}\nolimits}\ignorespaces}
\def\makeoplist#1 {\def\@@tmpa{#1}\def\@@tmpb{***}%
  \ifx\@@tmpa\@@tmpb\else\makeop{#1}\expandafter\makeoplist\fi}
\makeoplist 
ad Ad Alb Alt Aut Br can Card card Char Cl Coker coker Cond cond Corr
depth diag Diff Disc disc Div div dom End Ext Fil Fr Frac Frob Gal GL
Gr gr GSp GSpin Hom hom Id id Im im Ind ind Int inv Irr irr Isom
Ker ker len length Lie Nr Nrd O Ob ob ord PGL Pic Pr pr Proj PSL rank
RD Res Sgn sgn SL SO Sp Sp sp Spec Spf Spin 
St st Stab stab Std std SU Supp supp Swan Sym Tor tor tors 
torsion Tr tr Trd U val Vol vol ***
\def\makermlist#1 {\def\@@tmpa{#1}\def\@@tmpb{***}\ifx\@@tmpa\@@tmpb
  \else\expandafter\def\csname#1\endcsname{{\rm#1}}\expandafter\makermlist\fi}
\makermlist ab alg arith can geom opp op red sep sh ***
\makeatother
\DeclareMathAlphabet\eusm{U}{eus}{m}{n}
\def\makebb#1{\expandafter\def\csname bb#1\endcsname{{\mathbb{#1}}}\ignorespaces}
\def\makerm#1{\expandafter\def\csname rm#1\endcsname{{\rm #1}}\ignorespaces}
\def\makebf#1{\expandafter\def\csname bf#1\endcsname{{\bf #1}}\ignorespaces}
\def\makegr#1{\expandafter\def\csname gr#1\endcsname{{\mathfrak{#1}}}\ignorespaces}
\def\makescr#1{\expandafter\def\csname scr#1\endcsname{{\mathscr{#1}}}\ignorespaces}
\def\makecal#1{\expandafter\def\csname cal#1\endcsname{{\cal #1}}\ignorespaces}
\def\makeudl#1{\expandafter\def\csname udl#1\endcsname{{\underline{#1}}}\ignorespaces}
\def\doLetters#1{%
  #1A #1B #1C #1D #1E #1F #1G #1H #1I #1J #1K #1L #1M
  #1N #1O #1P #1Q #1R #1S #1T #1U #1V #1W #1X #1Y #1Z}
\def\doletters#1{%
  #1a #1b #1c #1d #1e #1f #1g #1h #1i #1j #1k #1l #1m
  #1n #1o #1p #1q #1r #1s #1t #1u #1v #1w #1x #1y #1z}
\doLetters\makebb\doLetters\makecal\doLetters\makerm\doLetters\makebf\doLetters\makescr
\doletters\makerm \doletters\makebf \doLetters\makegr \doletters\makegr \doLetters\makeudl \doletters\makeudl
\def\proof{\medbreak\noindent{\scshape Proof.}\enspace}%
     \def\qed{\qedmark\medbreak}%
\def\qedmark{{\enspace\vrule height 6pt width 5pt depth 1.5pt}}%
    \def\setminus{\smallsetminus}

\def\injto{\hookrightarrow}

\makeatletter   \newdimen\mina@@\mina@@=18pt
\newcommand{\xrtarw}[2][]{\mathrel{\mathop{\,\setbox\z@\vbox{\m@th
  \hbox{$\scriptstyle\;{#1}\;\;$}\hbox{$\m@th\scriptstyle\;{#2}\;\;$}}%
  \hbox to\ifdim\wd\z@>\mina@@\wd\z@\else\mina@@\fi{\rightarrowfill@
  \displaystyle}\,}\limits^{#2}\@ifnotempty{#1}{_{#1}}}}
\newcommand{\xltarw}[2][]{\mathrel{\mathop{\,\setbox\z@\vbox{\m@th
  \hbox{$\scriptstyle\;\;{#1}\;$}\hbox{$\m@th\scriptstyle\;\;{#2}\;\;$}}%
  \hbox to\ifdim\wd\z@>\mina@@\wd\z@\else\mina@@\fi{\leftarrowfill@
  \displaystyle}\,}\limits^{#2}\@ifnotempty{#1}{_{#1}}}}
\makeatother
\def\XYmatrix{\xymatrix@M=5pt} 
\def\ncmd{\newcommand}
\ncmd{\xysubset}[1][r]{\ar@<-2.5pt>@{^(-}[#1]\ar@<2.5pt>@{_(-}[#1]}
\ncmd{\XYmatrixc}[1]{\vcenter{\XYmatrix{#1}}}
\ncmd{\xyto}[1][r]{\ar@{->}[#1]}      \ncmd{\xyinj}[1][r]{\ar@{^(->}[#1]}
\ncmd{\xysurj}[1][r]{\ar@{->>}[#1]}   \ncmd{\xyline}[1][r]{\ar@{-}[#1]}
\ncmd{\xydotsto}[1][r]{\ar@{.>}[#1]}  \ncmd{\xydots}[1][r]{\ar@{.}[#1]}
\ncmd{\xyleadsto}[1][r]{\ar@{~>}[#1]} \ncmd{\xyeq}[1][r]{\ar@{=}[#1]}
\ncmd{\xyequal}[1][r]{\ar@{=}[#1]}    \ncmd{\xyequals}[1][r]{\ar@{=}[#1]}
\ncmd{\xymapsto}[1][r]{\ar@{|->}[#1]}\ncmd{\xyimplies}[1][r]{\ar@{=>}[#1]}
\ncmd{\xytofrom}[1][r]{\ar@{<->}[#1]} 
\def\xyinjto{\xyinj}   
\def\XYTOTO[#1]^#2_#3{\xyto[#1]<0.5ex>^{#2}\xyto[#1]<-0.5ex>_{#3}}
\ncmd{\xytoto}[1][r]{\XYTOTO[#1]}
\def\beginmat{\begin{pmatrix}}\def\endmat{\end{pmatrix}}
\def\pmat#1]{{\def\beginmat{\begin{pmatrix}}\def\endmat{\end{pmatrix}}\mat#1]}}
\def\bmat#1]{{\def\beginmat{\begin{bmatrix}}\def\endmat{\end{bmatrix}}\mat#1]}}
\def\Bmat#1]{{\def\beginmat{\begin{Bmatrix}}\def\endmat{\end{Bmatrix}}\mat#1]}}
\def\vmat#1]{{\def\beginmat{\begin{vmatrix}}\def\endmat{\end{vmatrix}}\mat#1]}}
\def\Vmat#1]{{\def\beginmat{\begin{Vmatrix}}\def\endmat{\end{Vmatrix}}\mat#1]}}
\def\smat#1]{{\def\beginmat{\begin{smallmatrix}}%
  \def\endmat{\end{smallmatrix}}\left(\mat#1]\right)}}
\def\mat#1#2]{\ifcase#1\or \matA#2]\or \matAA#2]\or \matAAA#2]\fi}
\def\matA  #1#2]{\ifcase#1\or \matAB  #2]\or \matABB  #2]\or \matABBB  #2]\fi}
\def\matAA #1#2]{\ifcase#1\or \matAAB #2]\or \matAABB #2]\or \matAABBB #2]\fi}
\def\matAAA#1#2]{\ifcase#1\or \matAAAB#2]\or \matAAABB#2]\or \matAAABBB#2]\fi}
\def\matAB[#1]{\beginmat#1\endmat}
\def\matABB[#1,#2]{\beginmat#1&#2\endmat}
\def\matABBB[#1,#2,#3]{\beginmat#1&#2&#3\endmat}
\def\matAAB[#1;#2]{\beginmat#1\\#2\endmat}
\def\matAABB[#1,#2;#3,#4]{\beginmat#1&#2\\#3&#4\endmat}
\def\matAABBB[#1,#2,#3;#4,#5,#6]{\beginmat
   #1&#2&#3\\#4&#5&#6\endmat}
\def\matAAAB[#1;#2;#3]{\beginmat#1\\#2\\#3\endmat}
\def\matAAABB[#1,#2;#3,#4;#5,#6]{\beginmat
   #1&#2\\#3&#4\\#5&#6\endmat}
\def\matAAABBB[#1,#2,#3;#4,#5,#6;#7,#8,#9]{\beginmat
   #1&#2&#3\\#4&#5&#6\\#7&#8&#9\endmat}
\def\beginalignorgather#1#2\endalignorgather{%
  \ifx#1!\beginaorgnostar#2\endaorgnostar\else\beginaorgstar#1#2\endaorgstar\fi}
\def\beginaorgstar#1#2\endaorgstar{%
  \ifx#1@\begin{align*}#2\end{align*}\else\begin{gather*}#1#2\end{gather*}\fi}
\def\beginaorgnostar#1#2\endaorgnostar{%
  \ifx#1@\begin{align}#2\end{align}\else\begin{gather}#1#2\end{gather}\fi}
\def\[#1\]{\beginalignorgather#1\endalignorgather}

\def\dbltag#1#2{\tag*{\hbox to 0pt{\hbox to \hsize
  {\hfil#2}\hss}#1}}

\newcommand{\lowsim}{\smash{\hbox{\lower2.5pt
  \hbox{\(\scriptstyle\sim\)}}}}
\makeatletter
\def\itemize{
  \ifnum \@itemdepth >\thr@@\@toodeep\else\advance\@itemdepth\@ne
  \edef\@itemitem{labelitem\romannumeral\the\@itemdepth}\expandafter\list
  \csname\@itemitem\endcsname{\def\makelabel##1{\hss\llap{##1}}}\fi\itemsep=-2pt}
\makeatother

\def\SEC#1#2\par{{\Large\bigbreak\noindent\llap{\bf #1 }{\bf #2}\par\medskip}}
\def\PARA#1. {\medbreak\noindent{\bf #1.}\enspace\ignorespaces}

\def\TH#1{\medbreak\noindent{\bfseries\scshape #1}\enspace
  \bgroup\it\ignorespaces}
\def\endTH{\egroup\medbreak}

\def\PAR{\medbreak\indent\llap{\hbox to \parindent
  {$\blacktriangleright$\hfil}}\ignorespaces}
\newcommand{\Hecke}{\mathscr{H}}
\newcommand{\rK}{\mathcal K}
\newcommand{\rKp}{\mathcal K_+}

\usepackage[all]{xy}

\begin{document}

\begin{center}\scshape
\LARGE
  Construction of tame types 
\medbreak
\normalsize
  Ju-Lee Kim\footnote{partially supported by an NSF grant}\\
\smallbreak
\footnotesize
  Massachusetts Institute of Technology\\
\smallbreak
  and 
\smallbreak
\normalsize
  Jiu-Kang Yu\\
\footnotesize
 Chinese University of Hong Kong
\end{center}
\bigbreak

\begin{center}
{\it To Professor Roger Howe on the occasion of his 70th Birthday}
\end{center}

\bigskip

\section{Introduction}

Let $F$ be a non archimedean local field and $G$ a reductive group defined over $F$. To study the category \(\scrR:=\scrR(G)\) of smooth representations of $G$, Bushnell and Kutzko \cite{BK} formulated the theory of types, which build on the theory of minimal $K$-types \cite{HM} and a general framework due to Bernstein \cite{Be}.  This
category $\scrR$ is decomposed into a  product of additive
subcategories by the theory of Bernstein center \cite{Be}:
\begin{equation}\label{i: cat dec}
\scrR=\prod_{\grs \in \scrI} \scrR^\grs(G),
\end{equation}
where $\scrI$ is the set of inertial classes of cuspidal pairs.
Moreover, each factor \(\scrR^\grs:=\scrR^\grs(G)\) can not be further decomposed (see \cite{BK} for details).  When we have a type in the sense of Bushnell and Kutzko, we have the means to study a finite number of these \(\scrR^\grs\).  To be
effective, we would like to have one type for each single
\(\scrR^\grs\).  The existence and the nature of such types has
been a fundamental problem.

In this article, we will give a construction of types for a general
reductive \(p\)-adic group \(G\).  
Our method produces types corresponding to a single factor
\(\scrR^\grs\).  The construction here is not new: it is the same
construction used in \cite{Yu} for supercuspidal representations, when
certain obvious constraints pertaining to supercuspidality are removed.
What is new is that an additional constraint must be imposed
concerning the embeddings of buildings.  For any (tamely ramified
twist of a) Levi subgroup \(G'\) of \(G\), there is a family of
embeddings of the (extended) building \(\scrB(G')\) of \(G'\) into that \(\scrB(G)\) of $G$.
This family forms a Euclidean space.  The choice of embeddings is
unimportant for almost all applications.  But it is crucial here
that we avoid a certain set of embeddings of measure \(0\). 
This generic choice of embeddings allows us to construct our types as covers of supercuspidal types on Levi subgroups in a uniform manner. To prove that this construction indeed yields $G$-covers and thus types in the sense of Bushnell and Kutzko \cite{BK}, we use ideas from the work of \cite{Kim, MP2}.

In \S\ref{sec: exhaustion}, we sketch a proof that our construction
yields sufficiently many types to study all irreducible admissible
representations under a suitable ``tameness'' hypothesis on \(G\) and $F$. 

\vspace{6pt}

\noindent{\bf Acknowledgment} It is our great pleasure to dedicate this paper to Professor Roger Howe.
We would like to thank the referee for his/her thoughtful and detailed comments, which have been invaluable in improving the exposition. An earlier version of this paper (without Sections 9 and 10) had been circulated among experts last 10 years. We would like to thank our friends and colleagues who encouraged us to publish the paper.

\vspace{6pt}

\section{Notation and conventions}

\subsection[] We adopt all notation and conventions from \cite[p.~582]{Yu}.
However, in this paper we do not need to treat base field
extensions extensively except in (\ref{P-iwa}).  Therefore, we work over a fixed
non-archimedean local field \(F\), that is, $F$ is either a $p$-adic field or a function field over finite field. If \(G\) is an algebraic group
over \(F\), we will denote $G$'s group of rational points also by \(G\) for simplicity.
This should lead to no confusion.

\subsection[] Throughout this paper, \(G\) is a connected reductive
group over \(F\), split over a tamely ramified extension of \(F\). For any maximal $F$-split torus $S$, $\Phi(G,S,F)$ denotes the corresponding set of roots in $G$. For $a\in\Phi(G,S,F)$, let $U_a$ (resp. $\mathfrak u_a$) be the root subgroup (resp. root subspace) corresponding to $a$.

By a Levi subgroup of \(G\), we mean an \(F\)-subgroup of \(G\) which
is a Levi factor of a parabolic \(F\)-subgroup of \(G\).  By a twisted
Levi subgroup \(G\), we mean an \(F\)-subgroup \(G'\) of \(G\) such
that \(G' \otimes_F \bar F\) is a Levi subgroup of \(G \otimes_F \bar
F\).

\subsection[] \label{torsion}
We assume that the residue characteristic \(p\) of
\(F\) is not a torsion prime for \(\psi(G)^\vee\), the root datum dual
to the root datum \(\psi(G)\) of \(G \otimes_F \bar F\).  See \cite[\S7]{Yu} and 
\cite{St} for the relevant notions.  By \cite[2.3]{St},
\(p\) is not a torsion prime for \(\psi(G')^\vee\), for any (twisted) Levi
subgroup \(G'\) of \(G\).  From \S\ref{S-type} on, we also
assume that \(p\) is odd.

\subsection[]\label{levi-seq}
Let \(\vec G=(G^0,\ldots,G^d)\) be a tamely ramified twisted Levi
sequence in \(G\), that is, each $G^i$ is a $E$-Levi subgroup of $G$ over a tamely ramified finite extension $E$ of $F$ (\cite[p.~586]{Yu}).
Let \(M^0\) be an $F$-Levi subgroup of \(G^0\).  Let \(\scrZ_\rms(M^0)\) be the
maximal \(F\)-split torus of the center \(\scrZ(M^0)\) of \(M^0\).  We
define \(M^i\) to be the centralizer of \(\scrZ_\rms(M^0)\) in \(G^i\).

\begin{lemma*} \ 
  \begin{itemize}
  \item [\rm(a)] \(M^i\) is an \(F\)-Levi subgroup of \(G^i\).
  \item [\rm(b)] \(\vec M:=(M^0, M^1, \ldots, M^d)\) is a generalized
    twisted Levi sequence in \(M:=M^d\) in the sense of [Yu, page 616].
  \item [\rm(c)] \(\scrZ(M^0)/\scrZ(M^d)\) is \(F\)-anisotropic.
  \end{itemize}
\end{lemma*}

\proof  (a) follows from \cite[20.4]{Bo}.  It then follows that \(M^i\)
is a twisted Levi subgroup of \(G^j\) for \(i \leq j\).  Therefore,
for \(i \leq j\), \(M^i\) is the centralizer of \(\scrZ(M^i)^\circ\) in
\(G^j\), hence is also the centralizer of \(\scrZ(M^i)^\circ\) in
\(M^j\).  Again by \cite[20.4]{Bo}, this implies that \(M^i\) is a
twisted Levi subgroup of \(M^j\).  We have proved (b).

Finally, since \(\scrZ(M^d) \subset \scrZ(M^0)\), the \(F\)-split rank of
\(\scrZ(M^d)^\circ\) is smaller than or equal to that of
\(\scrZ(M^0)^\circ\).  By construction,
\(\scrZ(M^d)^\circ\supset\scrZ_\rms(M^0)\).  Therefore
\(\scrZ(M^d)^\circ\) and \(\scrZ(M^0)^\circ\) have the same \(F\)-split
rank.  This proves (c).\qed

\begin{rem*}
We observe that $\vec M$ is a tamely ramified twisted Levi sequence in $M$.
\end{rem*}

\section{Generic embeddings of buildings}

\subsection[] We recall that, if \(G'\) is a tamely ramified twisted Levi subgroup of \(G\), then
there exists a family of natural embeddings of buildings \(\scrB(G') \injto \scrB(G)\),
which is an affine space under \(X_*(\scrZ_\rms(G')) \otimes \bbR\).
All these embeddings have the same image.  
Two embeddings in the
same orbit of \(X_*(\scrZ_\rms(G)) \otimes\bbR\) can be regarded as
the same for most purposes.

\begin{defn}
  Let \(M\) be a Levi subgroup of \(G\), \(y \in \scrB(M)\), and \(s
  \in \bbR\).  We say that the embedding \(\iota:\scrB(M) \injto
  \scrB(G)\) is {\it \((y,s)\)-generic}, or {\it \(s\)-generic with
    respect to \(y\)}, if \(U_{a,\iota(y),s}=U_{a,\iota(y),s+}\) for
  all \(a \in \Phi(G,S,F) \setminus \Phi(M,S,F)\), where \(S\) is any
maximal \(F\)-split torus of \(M\) such that \(y\in A(M,S,F)\). Here $A(M,S,F)$ is the apartment associated to $S$ in $\scrB(M)$.

Once an embedding $\iota$ is fixed, we will identify $\scrB(M)$ as a subset of $\scrB(G)$.
\end{defn}

Here, \(\{U_{a,\iota(y),r}\}_{r \in \bbR}\) is the filtration on the root group \(U_a\), $a\in\Phi(G,S,F)$ so that 
\(U_{a,\iota(y),r}=U_a \cap G_{\iota(y),r}\), where
\(\{G_{\iota(y),r}\}_{r \geq 0}\) is the Moy-Prasad filtration 
(see \cite{MP1}, \cite{MP2}).  The
following two results illustrate the usefulness of the notion of generic
embeddings.

\begin{prop} \label{P-mp2}
  Let \(G,M,y\) be as above and let \(\iota:\scrB(M)
  \injto \scrB(G)\) be \(0\)-generic relative to \(y\).  Let \(P=MU\)
  be a parabolic \(F\)-subgroup of \(G\) with Levi factor \(M\).
  For any smooth representation \(V\) of \(G\), the natural map \(r_U:
  V \to V_U\) from \(V\) to its Jacquet module induces a bijection
\[
r_U: V^{G_{y,0+}} \to (V_U)^{M_{y,0+}}.
\]
\end{prop}

This is a reformulation of \cite[Proposition 6.7]{MP2}. Note that $y\in\scrB(M)$ and an embedding $\scrB(M)\hookrightarrow\scrB(G)$ is not always $0$-generic with respect to $y$.

\begin{rem} \label{rmk: vdash} We can use generic embeddings to gain some new insight 
for the result in \cite[\S17]{Yu}.  Indeed, when $G'$ is a Levi subgroup, one observes that the main
result in \cite[Theorem 17.1]{Yu} is obvious when the embedding \(\iota :\scrB(G')
\injto \scrB(G)\) implicitly used is \((y,s)\)-generic, where \(s=r/2\).
In this case, one can argue directly using the last paragraph of the proof \cite[Corollary 17.3]{Yu}.

In general, one argues that there is an embedding \(\iota_1:\scrB(G')
\injto \scrB(G)\) close to \(\iota\) which is \((y,s)\)-generic (see
(\ref{L-gen}) below), and \((J,\tilde
\phi)=(J,\ind_{J_1}^J\tilde\phi_1)\) where \((J_1,\tilde\phi_1)\) is
constructed in the same way as \((J,\tilde\phi)\) but using
\(\iota_1\) in place of \(\iota\).  Then the theorem follows
immediately from the generic case.  In fact, this is just rephrasing
the proof in \cite[\S17]{Yu}.  Our \(J_1\) is the ad hoc object
\(J_\vdash\) used there. However, now we view \cite[Lemma 17.2]{Yu} as a
literal special case of \cite[Theorem 9.4]{Yu} by varying the embedding,
and we should regard \cite[Theorem 17.1]{Yu} in the case of a generic
embedding as the essential result.
\end{rem}


\subsection[] \label{emb}
We now work in the setting of \ref{levi-seq}.  Consider
a commutative diagram of embeddings:
\[\tag*{\(\{\iota\}:\)}\vcenter{
\XYmatrix{
\scrB(G^0) \xyinjto[r] &\scrB(G^1) \xyinjto[r] &\cdots\xyinjto[r]
&\scrB(G^d)\\
\scrB(M^0) \xyinjto[r]\xyinjto[u] &\scrB(M^1) \xyinjto[r]\xyinjto[u]
 &\cdots\xyinjto[r] &\scrB(M^d) \xyinjto[u]
}}
\]
To specify such a diagram of embeddings,
it suffices to give the image of a fixed
\(y \in \scrB(M^0)\) in \(\scrB(M^i)\), \(1 \leq i \leq d\) and
in \(\scrB(G^i)\), \(0 \leq i \leq d\).  We will denote this
whole diagram by \(\{\iota\}\), and we will denote by \(\iota\) any
composite embedding in this diagram (from \(\scrB(M^i)\) to
\(\scrB(M^j)\) or \(\scrB(G^j)\), or from \(\scrB(G^i)\) to
\(\scrB(G^j)\), for \(i \leq j\)).

\begin{defn*}  Let \(\vec s=(s_0,\ldots,s_{d})\) be a sequence of real
  numbers, and \(y \in \scrB(M^0)\). We say that \(\{\iota\}\) is \(\vec
  s\)-{\it generic (relative to \(y\))} if \(\iota:\scrB(M^i) \to \scrB(G^i)\) is
  \(s_{i}\)-generic relative to \(\iota(y) \in \scrB(M^i)\) for \(0 \leq
  i \leq d\).
\end{defn*}

\subsection[]\label{L-gen}
We now establish the abundance of
generic embeddings.  

Let \(\{\iota\}\) be a commutative diagram of embeddings
as in \S\ref{emb}, and \(y \in \scrB(M^0)\).  Denote the image of \(y\) in
\(\scrB(M^i)\) by \(y_i\), and that in
\(\scrB(G^i)\) by \(z_i\), \(0 \leq i \leq d\).  Let \(v \in
X_*(\scrZ_\rms(M^0)) \otimes\bbR\).  There is a commutative diagram of
embeddings, to be denoted by \(\{\iota\}_v\), in which the image of
\(y\) in \(\scrB(M^i)\) is \(y_i\), and that in \(\scrB(G^i)\) is
\(z_i+v\), \(0 \leq i \leq d\).

\begin{lemma*}\label{lem: genemb} Fix \(\vec s \in \bbR^{d+1}\).
\begin{itemize}
\item[\rm (a)] \(\vec s\)-generic commutative diagrams of embeddings exist.
\item[\rm (b)] Assume that \(G \neq M\). For \(0 \leq i \leq d\), let $S^i$ be a maximal \(F\)-split torus of \(M^i\).  Let \(\gamma \in X_*(\scrZ_\rms(M^0)) \otimes \bbR\) be such
  that \(\langle a,\gamma\rangle \neq0\) for \(a
  \in \Phi(G^i,S^i,F) \setminus \Phi(M^i,S^i,F)\), \(0 \leq i \leq d\).
  Then the set of \(t \in \bbR\) such that the commutative
  diagram of embeddings \(\{\iota\}_{t\gamma}\) is not \(\vec s\)-generic
  is an infinite discrete subset of \(\bbR\).
\end{itemize}
\end{lemma*}

\proof 
For \(i=0,1,\cdots,d\) and \(a\in\bigl(\Phi(G^i,S^i,F) \setminus
\Phi(M^i,S^i,F)\bigr)\), there exist infinite discrete subsets
\(\Gamma _{i,a}\) of \(\bbR\) such that the set of \(v \in V=X_*(\scrZ_\rms(M^0))\otimes\bbR\)
with \(\{\iota\}_v\) not \(\vec s\)-generic is the union of
hyperplanes in \(V\)
defined by \(a(v)=c\), \(c \in \Gamma_{i,a}\).
Both statements follow easily from this.\qed

\section{Covers and decompositions}

\subsection[Iwahori-type decompositions] Let \(P=MU\) be a parabolic
\(F\)-subgroup of $G$ with Levi factor \(M\), and \(\bar P=M\bar U\) the
opposite parabolic.  A compact open subgroup \(K\) of \(G\) is said to
{\it decompose with respect to\/} \(U,M,\bar U\) if 
\[
K=(K \cap U).(K \cap M).(K \cap \bar U).
\]

\subsection[Covers] \label{D-cover}
Let \(M\) be a Levi subgroup of \(G\), \(K\) (resp.~\(K_M\)) a compact
open subgroup of \(G\) (resp.~\(M\)), and \(\rho\) (resp.~\(\rho_M\))
an irreducible smooth representation of \(K\) (resp.~\(K_M\)).  The pair
\((K,\rho)\) is called a {\it
  \(G\)-cover} of \((K_M,\rho_M)\) if for any opposite pair of
parabolic subgroups \(P=MU, \bar P=M\bar U\) with Levi factor \(M\),
we have
\begin{itemize}
\item [(i)] \(K\) decomposes with respect to \((U,M,\bar U)\).
\item [(ii)] \(\rho|K_M=\rho_M\) and \(K \cap U, K \cap \bar U \subset
  \ker(\rho)\).
\item [(iii)] For any smooth representation \(V\) of \(G\), the natural
  map from \(V\) to its Jacquet module \(V_U\) induces an injection on
  \(V^{(K,\rho)}\), the \((K,\rho)\)-isotypic subspace of \(V\).
\end{itemize}
This definition is due to Bushnell and Kutzko \cite{BK}, although we have
used a reformulation given in \cite[Th\'eor\`eme 1]{Bl} (see also \cite[\S4.1]{GR}).

\subsection\label{P-iwa}
We now give a useful (probably well-known) class of compact open
subgroups with the decomposition property with respect to \((U,M,\bar U)\).

Fix \(\iota:\scrB(M) \injto \scrB(G)\) and consider \(\scrB(M)\) as a
subset of \(\scrB(G)\).  Let \(E/F\) be a finite Galois extension such
that \(G \otimes E\) is split.  Let \(T\) be a maximal torus of \(M\),
defined over \(F\) and split over \(E\).  Let \(y\) be a point of
\(A(M,T,E) \cap \scrB(G)\).  Put \(\Phi^0=\Phi(G,T,E) \cup \{0\}\) and
let \(f:\Phi^0 \to \tilde \bbR\) be a \(\Gal(E/F)\)-stable and concave
function.  Then we can define \(G(E)_{y,f}\) and
\(K=G_{y,f}=G(E)_{y,f} \cap G\) as in \cite[page 608]{Yu}.

In addition, let \(K_+=G_{y,f_+}\), where \(f_+\) is the concave
function \(a \mapsto \max(f(a),0+)\) \cite[Lemma 13.1 (ii)]{Yu}.

\begin{prop*} 
  Suppose that \(f(a) \geq 0\) for all \(a \in \Phi^0\)
  and \(f(a)>0\) for all \(a \in \Phi^0\setminus
  \Phi(M,T,E)\).
  \begin{itemize}
  \item [\rm (a)] \(K=G_{y,f}\) and \(K_+=G_{y,f_+}\) decompose with respect to
  \((U,M,\bar U)\).
  \item [\rm (b)] Let \(\hat K_M\) be a compact open subgroup of \(M\) containing
  \(K_+ \cap M\) such that \(\hat
  K_M\) normalizes \(K_+\).  Then \(\hat K:=(K \cap U)\hat K_M(K\cap
  \bar U)\) is an open compact subgroup, which decomposes with respect
  to \((U,M,\bar U)\).  Moreover, \(\hat K/K_+ \simeq \hat K_M/(K_+
  \cap M)\).
  \end{itemize}
\end{prop*}

\proof (a) It suffices to prove the assertion for \(K_E=G(E)_{y,f}\).
Indeed, if \(K_E\) decomposes with respect to \((U,M,\bar U)\),
and \(g \in G_{y,f}\), then \(g = u m \bar u\) for unique elements \(u
\in K_E \cap U(E), m \in K_E \cap M(E), \bar u \in K_E \cap \bar
U(E)\).  For any \(\sigma  \in \Gal(E/F)\), \(um\bar u=g=\sigma(g)=\sigma
(u)\sigma (m)\sigma (\bar u)\).  By the uniqueness of the
decomposition (a consequence of the big cell theorem \cite[14.21]{Bo}),
we have \(u \in U(E)^{\Gal(E/F)}=U(F)\), etc.  It follows that
\(K=G_{y,f}\) decomposes with respect to \((U,M,\bar U)\).  A similar
remark applies to \(K_+\).

We now prove (a) with the additional assumption
that \(E=F\) and \(T\) is split over \(F\).  The statement about
\(K_+\) is then a special case of \cite[6.4.48]{BT1}.
Write \(\Phi^0=\Phi_U
\sqcup \Phi_M \sqcup \Phi_{\bar U}\), where \(\Phi_M=\Phi(M,T,F)
\cup \{0\}\), and \(\Phi_U\) (resp.~\(\Phi_{\bar U}\)) is the set of roots
for the action of \(T\) on \(\Lie U\) (resp.~on \(\Lie \bar U)\)).
For \(H=U,M,\bar U\), let 
\[
f_H(a)=\begin{cases}
f(a) &\mbox{ if \(a \in \Phi_H\),}\\
\infty&\mbox{otherwise.}
\end{cases}
\]
Then \(f_H\) is concave by \cite[Lemma 13.1 (iv)]{Yu}, and 
\(K^H=G_{y,f_H} \subset K \cap H\).

By \cite[6.4.43]{BT1}, \(K^M\) normalizes \(K^U\).  Therefore, \(K^U K^M\)
is a subgroup of \(G\), and is the same as \(G_{y,f_{UM}}\), where
\(f_{UM}=\inf(f_U,f_M)\).  Again by \cite[6.4.43]{BT1}, \(K^U K^M \subset
K\) normalizes \(K_+\).  Therefore, \(K^U K^M K_+\) is a subgroup of
\(G\).  Since \(K_+ = K^U (K_+ \cap M) K^{\bar U}\) and \(K^U (K_+
\cap M) \subset K^U K^M\), we have \(K^U K^M K_+ = K^U K^M K^{\bar
  U}\) is a subgroup of \(G\).  It follows that this subgroup is
\(G_{y,f}=K\), and \(K^H=K \cap H\) for \(H=U,M,\bar U\).  This proves
(a).

(b) Since \(\hat K_M\subset M\) normalizes \(U\), \(\hat K_M\)
normalizes \(K^U\) and \(K^U\hat K_M\) is a subgroup of \(G\).  This
subgroup normalizes \(K_+\) since both \(K^U\) and \(\hat K_M\) do.
Therefore, \(K^U \hat K_M K_+\) is a subgroup of \(G\).  Clearly, this subgroup
is equal to \(K^U \hat K_M K^{\bar U}\).

The natural morphism \(\hat K_M \to \hat K/K_+\) is surjective with
kernel \(\hat K_M \cap K_+=(K_+ \cap M)\).  This gives the asserted
isomorphism and finishes the proof of (b).\qed

\section{Heisenberg triples}

\begin{defn} \label{def-h3}
  Let \(J \supset J_+\) be compact, open, pro-\(p\) subgroups of
  \(G\), and let \(\varphi:J_+ \to \bbC^\times\) be a smooth character
  such that \(\varphi(J_+)=\mu_p:=\{\zeta\in\bbC^\times:\zeta^p=1\}\).
  We say that \((J,J_+,\varphi)\) is a Heisenberg triple if
  \begin{itemize}
  \item [(i)] \(J_+\) is a normal subgroup of \(J\) and 
    \(J/J_+\) is an abelian group of exponent \(p\).
  \item [(ii)] The commutator subgroup \([J,J_+] \subset \ker(\varphi)\).
  \item [(iii)] The symplectic pairing \((J/J_+) \times (J/J_+) \to
    \mu_p\), \((a,b) \mapsto \varphi(aba^{-1}b^{-1})\) is non-degenerate.
  \end{itemize}
\end{defn}

Notice that (i) and (ii) imply that the pairing in (iii) is well-defined.  
It follows that \(J/\ker(\varphi)\) is a Heisenberg \(p\)-group.  Such
triples often occur in the representation theory of \(p\)-adic groups.

\begin{exam} \label{h-old}
  We now recall the fundamental Heisenberg triple used in \cite{Yu}.  The
  setting of this example will be in force through the rest of this section.
  Let \((G',G)\) be a tamely ramified twisted Levi sequence in \(G\), and
  \(y \in \scrB(G')\).  Fix an embedding \(\scrB(G') \injto \scrB(G)\)
  to identify \(\scrB(G')\) as a subset of \(\scrB(G)\).
  Let \(r\) be a positive real number and \(\phi:G'_{y,r:r+}\to\bbC^\times\) a
  \(G\)-generic character in the sense of \cite[\S9]{Yu}.  Put \(s=r/2\),
  \(J=(G',G)_{y,(r,s)}\), \(J_+=(G',G)_{y,(r,s+)}\),
  and \(\varphi:J_+\to\bbC^\times\) the character obtained by extending the restriction of \(\phi\) to $G'_{y,r}$ trivially across  a subgroup of $G$ that is `perpendicular' to $G'$ in a suitable sense (see \cite[\S4]{Yu}).  Then \((J,J_+,\varphi)\) is a Heisenberg
  triple by \cite[Lemma 11.1]{Yu}.
\end{exam}

\subsection\label{subsec: emb} Keep the settings in (\ref{h-old}).  In addition,
let \(M'\) be a Levi subgroup of \(G'\) such that $y\in\scrB(M')$. Let $M$ be the centralizer in $G$ of $\scrZ_s(M')$.  Then, if we put \((G^0,G^1)=(G',G)\), $(M^0,M^1):=(M',M)$ is a Levi sequence as in \ref{levi-seq}.  We also fix a commutative diagram of embeddings of buildings extending $\scrB(G') \injto \scrB(G)$:
\[\XYmatrix{\scrB(G') \xyinjto[r] &\scrB(G)\\
\scrB(M') \xyinjto[r]\xyinjto[u]&\scrB(M)\xyinjto[u]
}
\]
We treat these embeddings as inclusions.

\begin{lemma*} \(\phi_M:=\phi|M'_{y,r}\) is \(M\)-generic of depth \(r\)
  relative to \(y\).
\end{lemma*}

\proof By definition, \(\phi|G'_{y,r:r+}\) is realized by a
\(G\)-generic element \(X^* \in (\Lie^*\scrZ(G'))_{-r}\).  Let
\(E/F\) be a finite extension over which \(M',M,G'\) and \(G\) are
all split, and \(T\) a maximal \(E\)-split torus of \(M'\) such that
\(y \in A(M',T,E)\).  Then the genericity of \(X^*\) means:
\(\ord(X^*(H_a))=-r\) for all \(a \in \Phi(G,T,E) \setminus
\Phi(G',T,E)\) (recall that \(H_a=da^\vee(1)\) and \(a^\vee:\bbG_\rmm
\to T\) is the coroot of \(a\)).  This is condition {\bf GE1} in [Yu, \S8].
By (\ref{torsion}) and [Yu, Lemma 8.1], {\bf GE2} holds automatically.

Clearly, \(\phi_M|M'_{y,r:r+}\) is realized by \(X^*_M:=\) the image
of \(X^*\) under \(\Lie^* G' \to \Lie^* M'\).  So we have
\(X^*_M(H_a)=X^*(H_a)\) has valuation \(-r\) for all \(a \in \Phi(M,T,E) \setminus
\Phi(M',T,E) \subset \Phi(G,T,E)\setminus \Phi(G',T,E)\) (notice that,
for \(a \in \Phi(M,T,E)\subset \Phi(G,T,E)\), \(H_a\) is the same
whether we consider \(a\) as a root of \(M\) or of \(G\)).  Again by
(\ref{torsion}) and [Yu, Lemma 8.1], this shows that \(\phi_M\) is
\(M\)-generic of depth \(r\) relative to \(y\).\qed

\subsection \label{h-new} 
Keep the settings in (\ref{h-old}) and (\ref{subsec: emb}). We now give the crucial construction of a new
Heisenberg triple needed in this paper.  Let
\[
\tilde J=(M,G)_{y,(r,s)} \cap (G',G)_{y,(r,s)}, \qquad
\tilde J_{+}=(M,G)_{y,(r,s)} \cap (G',G)_{y,(r,s+)}, \qquad
\tilde \varphi=\varphi|\tilde J_+,
\]
where \(\varphi\) is the character used in \ref{h-old}.  By [Yu, Lemma
13.2], the groups \(\tilde J\) and \(\tilde J_+\) can be described by
concave functions in a suitable way.

\begin{lemma*} \ 
  \begin{itemize}
  \item [\rm (a)] \((\tilde J,\tilde J_+,\tilde\varphi)\) is a
    Heisenberg triple.
  \item [\rm (b)] Let \(P=MU\) be a parabolic \(F\)-subgroup of \(G\) with Levi factor
\(M\) and unipotent radical \(U\).  Let \(\bar U\) be the unipotent
radical of the opposite parabolic subgroup. Let \(W\) and \(\bar W\)
be the images of \(\tilde J^U:=\tilde J\cap U\) and
\(\tilde J^{\bar U}:=\tilde J \cap \bar U\) in \(\tilde V:=\tilde
J/\tilde J_+\) respectively.  Then \(W\) and \(\bar W\) are maximal
isotropic subspaces of \(\tilde V\), and \(\tilde V=W+\bar W\).  In
other words, \(W\) and \(\bar W\) form a complete polarization of
\(\tilde V\).
  \item[\rm (c)] Let \(C\) be the center of \(H:=\tilde J/\ker(\tilde
    \varphi)\). Let \(A,\bar A\) be the images of \(\tilde J^U\) and
    \(\tilde J^{\bar U}\) in \(H\).  Then \(A\cap C=\bar A\cap
    C=\{1\}\).
  \item [\rm (d)] \((G',G)_{y,(0,s)} \cap U\) normalizes \(\tilde J^U\).
  \end{itemize}
\end{lemma*}

\proof Let \((J,J_+,\varphi)\) be the triple in \ref{h-old}, so that
\(\tilde J \subset J\) and \(\tilde J_+ = \tilde J \cap
J_+\).  Thus conditions (i) and (ii) in
\ref{def-h3} follow from the corresponding statements for
\((J,J_+,\varphi)\).  It is also clear that \(\tilde\varphi(\tilde J_+) = \mu_p\).

Denote by \((J^M, J^M_+,\varphi^M)\) the Heisenberg triple
obtained by applying the construction of \ref{h-old} to \((M',M,y,\phi_M)\).
Then \(J^M \subset J\), \(J^M_+=J^M \cap
J_+\).  Put \(V=J/J_+\), \(V_M=J^M/J^M_+\),
\(\tilde V=\tilde J/\tilde J_+\).  Then we have embeddings of \(V_M\injto V\)
 and \(\tilde V\injto V\), which are compatible with the symplectic
pairings on these spaces.  We will regard these embeddings as inclusions.  

By [Yu, Lemma 13.3], \(J=J^M \tilde J\).  This implies \(V_M + \tilde
V=V\).  We now claim \(V_M \perp \tilde V\).  Since \(V_M\) is non-degenerate,
this will imply \(V = V_M \oplus \tilde V\) (orthogonal direct sum of
symplectic spaces) and \(\tilde V\) is non-degenerate.  So condition
(iii) in \ref{def-h3} will follow.

By Proposition~\ref{P-iwa},
\[
\tilde J=\tilde J^U M_{y,r} \tilde J^{\bar U}.
\]
It follows that we have \(\tilde V=W+\bar W\).  To
prove the claim, it
suffices to show that both \(W\) and \(\bar W\) are perpendicular to
\(V_M\) in \(V\).

Indeed, if \(a \in J^M\), \(b \in \tilde J^U\), then \(aba^{-1}b^{-1}
\in J_+\cap\tilde J^U\subset \tilde J_+ \cap \tilde J^U \subset \ker(\varphi)\).  This shows \(V_M \perp W\).
The same argument proves \(V_M \perp \bar W\).  This finishes the
proof of (a).

Similarly, if \(a, b \in \tilde J^U\), then \(aba^{-1}b^{-1} \in
\tilde J_+\cap \tilde J^U
\subset \ker(\varphi)\).  Therefore \(W\) is a totally isotropic
subspace of \(\tilde V\).  The same goes for \(\bar W\).  Since
\(V=W+\bar W\), both \(W\) and \(\bar W\) are maximal isotropic
subspaces.  This proves (b).

(c) is obvious since \(\tilde J_+ \cap \tilde J^U, \tilde
J_+\cap \tilde J^{\bar U} \subset \ker(\varphi)\).

To prove (d), it suffices to prove the analogous statement when \(F\)
is replaced by a finite, tamely ramified, Galois extension field
\(E\).  Therefore, we may and do assume that \(M'\) is split over
\(F\), and that there is a maximal \(F\)-split torus \(T\) of \(M'\) such
that \(y \in A(M',T,F)\).  Then we can write \(\tilde J^U=G_{y,f}\),
\((G',G)_{y,(0,s)} \cap U=G_{y,g}\) for suitable concave functions \(f,g\) on
\(\Phi^0=\Phi(G,T,F) \cup \{0\}\).

Explicitly, let \(v \in X_*(\scrZ_\rms(M)) \otimes \bbR\) be such that
\(\langle a,v\rangle >0\) for all roots \(a\) of \(\scrZ_\rms(M)\) on
the Lie algebra of \(U\).  Then 
\[
f(a)=\begin{cases}
r&\mbox{if \(\langle a,v\rangle >0,a \in \Phi(G',T,F)\),}\\
s&\mbox{if \(\langle a,v\rangle >0, a \notin \Phi(G',T,F)\),}\\
\infty&\mbox{if \(\langle a,v\rangle \leq0\);}
\end{cases}\qquad
g(a)=\begin{cases}
0&\mbox{if \(\langle a,v\rangle >0,a \in \Phi(G',T,F)\),}\\
s&\mbox{if \(\langle a,v\rangle >0, a \notin \Phi(G',T,F)\),}\\
\infty&\mbox{if \(\langle a,v\rangle \leq0\);}
\end{cases}
\]
According to [BT1, 6.4.43], it suffices to check
\[
f(pa+qb) \leq pf(a) + qg(b)
\]
whenever \(p,q\in\bbZ_{>0}\), \(a,b,pa+qb \in \Phi^0\).  This
condition is easily verified and hence (d) is proved.\qed

\begin{lemma} \label{L-nonzero}
  Let \(H\) be a finite Heisenberg \(p\)-group with center
  \(C\).  Assume that \(A\), \(\bar A\) are subgroups of \(H\) such
  that \(A \cap C=\bar A\cap C=\{1\}\), and the image of \(A\) and
  \(\bar A\) in \(V:=H/C\) form a complete polarization.  Let \(\psi\)
  be a non-trivial character of \(C\) and \((X,\rho)\) a complex
  representation of \(H\) such that \(\rho|C\) is
  \(\psi\)-isotypic.  Let \(v \in X^A\) be non-zero.  Then
  \(\sum_{b \in \bar A}b.v\) is also non-zero.
\end{lemma}

\proof This is \cite[Lemma 16.4]{Kim}.  For completeness, we produce a
proof here.  Assume \(v \neq 0\) is fixed by \(A\).  For \(b \in \bar A\), \(a \in
A\), we have
\[
ab.v=aba^{-1}b^{-1}ba.v=\psi(aba^{-1}b^{-1})(b.v).
\]
Therefore, \(b.v\) is an eigenvector for \(A\) for the character
\(\psi_b: a \mapsto \psi(aba^{-1}b^{-1})\).  As these characters are
distinct, the list of (non-zero) vectors \(\{b.v\}_{b \in \bar A}\) is linearly
independent.  It follows that their sum is non-zero.\qed

\section{Some covers of linear characters}

\subsection[Setup] \label{linear-setting}
We now work in the setting of \ref{levi-seq}.
Assume in addition: 
\begin{itemize}
\item[(i)] for each \(0 \leq i
\leq d-1\), we have a quasi-character \(\phi_i: G^i \to \bbC^\times\)
such that \(\phi_i\) is \(G^{i+1}\)-generic of depth \(r_i\) relative to any $x\in\scrB(G')$; and
\item[(ii)] these depths satisfy
\[
0 < r_0 < r_1 < \cdots <r_{d-1}; \mbox{ and}
\]
\item[(iii)] we have a point \(y \in \scrB(M^0)\) and
a commutative diagram of embeddings \(\{\iota\}\) as in \ref{emb},
which is \(\vec s\)-generic relative to \(y\), where
\[
\vec s=(0,s_0,\ldots,s_{d-1})=\left(0,\frac{r_0}{2},\frac{r_1}{2},\ldots,\frac{r_{d-1}}{2}\right)
\]
(notice that \(s_i=r_i/2\) is the \((i+1)\)-st component of \(\vec s\)
while \(r_i\) is the \(i\)-th component of \(\vec r\)).
\end{itemize}

We can now form compact subgroups similarly as in \cite[\S3]{Yu}:
\[
\rK^i
=G^0_{y,0}G^1_{y,s_0}\cdots G^i_{y,s_{i-1}},\qquad
\rKp^i=G^0_{y,0+}G^1_{y,s_0+}
\cdots G^i_{y,s_{i-1}+}.
\]

Let \(\theta_i:\rKp^i\to\bbC^\times\) be the character defined as in [Yu, \S4].
Put \(\rK=\rK^d\), \(\rKp=\rKp^d\), \(\theta=\theta_d\), \(\rKp^M=\rKp \cap M\),
\(\theta_M=\theta|\rKp^M\).  We caution the reader that
these groups do depend on the choice of \(\{\iota\}\), although the
dependency is suppressed in the above notation following [Yu].  We
will write \(\rK^i\{\iota\}\), \(\rKp^i\{\iota\}\), etc., when we need to make the dependency
clear.

\begin{rem*}

\begin{enumerate}
\item
In (i), if $\phi_i$ is $G^{i+1}$-generic of depth $r_i$ relative to an $y\in\scrB(G^i)$, it is $G^{i+1}$-generic of same depth relative to any $y\scrB(G^i)$. This follows from the definition of the genericity of $\phi_i$.
\item
Note that this \(\rK^i\) and $K^i$ in \S\ref{subsec: construction} are different from the \(K^i=G^0_{[y]}G^1_{y,s_0}\cdots G^i_{y,s_{i-1}}\) in \cite[p591]{Yu} in general. While $\rK^i$ and $K^i$ here are compact, $K^i$ in \cite{Yu} is compact mod center in general. 
\end{enumerate}

\end{rem*}

\begin{lemma}\label{tri-on-nil}
  Let \(P=MU\) be a parabolic \(F\)-subgroup of \(G\) with Levi factor
  \(M\) and unipotent radical \(U\), and \(\bar P=M\bar U\) the
  opposite parabolic.  Then we have
  \begin{itemize}
  \item [\rm (a)] \(\rK\) and \(\rKp\) decompose with respect to
    \((U,M,\bar U)\).
  \item [\rm (b)] The character \(\theta\) is trivial on \(\rK_+ \cap U\) and \(\rK_+ \cap \bar U\).
  \end{itemize}
\end{lemma}

\proof (a) follows from Proposition~\ref{P-iwa} (a).  We first
observe that \(\phi_i(U \cap G^i)=\{1\}\) since \(U \cap G^i\) lies in
the commutator subgroup of \(G^i\).  Assertion (b) follows from
this and the definition of \(\hat\phi_i\) and \(\theta_i\) in [Yu, \S4].\qed

\begin{thm}
  The pair \((\rKp,\theta)\) is a \(G\)-cover of \((\rK^M_+,\theta_M)\).
\end{thm}

\proof  We proceed by induction on \(d\).  The argument here is similar
to that in \cite[Proposition 17.2]{Kim}: 
the case \(d=0\) is a
reformulation of Proposition~\ref{P-mp2} and the inductive step
when \(d \geq 1\) follows the 
method of [MP2], using Lemma~\ref{L-nonzero} to
play the role of [MP2, Proposition 6.1].  

Now assume \(d \geq 1\) and also \(G \neq M\), since there is nothing
to prove if \(G=M\).  Among the three defining conditions for a cover
(\ref{D-cover}), (i) and (ii) are just the
preceding lemma.  It remains
to prove (iii): for any parabolic \(F\)-subgroup \(P=MU\) of \(G\) with Levi
factor \(M\), and for any smooth representation \(V\) of \(G\), 
the canonical map from \(V\) to its Jacquet module \(V_U\) is
injective on the \((\rKp,\theta)\)-isotypic subspace of \(V\).

Let \(v \in V\) be a non-zero \((\rKp,\theta)\)-isotypic vector.  It
suffices to show that \(\int_{N_i} g.v \, dg \neq 0\) for an increasing family of
open compact subgroups \(\{N_i\}\) of \(U\) whose union is the whole of \(U\).

Choose \(\gamma\in X_*(\scrZ_\rms(M^0)) \otimes\bbR\) such that
\(\langle a,\gamma\rangle >0\) for all roots \(a\) of \(\scrZ_\rms(M^0)\)
on the Lie algebra of \(U\).  For \(t \in \bbR\), form the groups
\(\rK(t)=\rK\{\iota\}_{t\gamma }\) and \(\rKp(t)=\rKp\{\iota\}_{t\gamma }\)
(see \ref{L-gen} for the notation).  Put \(N(t)=\rKp(t) \cap U\) and
\(\bar N(t)=\rKp(t) \cap \bar U\).  Then
\(N(t) \subset N(t')\), \(\bar N(t) \supset \bar N(t')\) for \(t < t'\) and \(\bigcup_{t \in \bbR} N(t)=U\).

By Lemma~\ref{L-gen}, there is an infinite sequence \(\cdots <t_{-1} <
t_0 < t_1 <
t_2 < \cdots\) such that \(t_n \to +\infty\) and \(t_{-n} \to -\infty\)
as \(n \to +\infty\), and \(\rKp(t)\) is constant
on the open intervals \(t_{i-1} < t < t_i\) (\(i \in \bbZ\)).
Therefore, we will denote by
\(\rKp(t_{i-1},t_i)\) the group \(\rKp(t)\) for any \(t \in (t_{i-1},t_i)\).

In fact, \(M\cap \rKp(t)=\rKp^M\), \(M \cap \rK(t)=\rK^M:=\rK\cap M\) for all
\(t\).  For \(i \in \bbZ\), let \(N_i=U \cap \rKp(t_{i-1},t_i)\), \(\bar N_i=\bar U
\cap \rKp(t_{i-1},t_i)\).  Then by Proposition~\ref{P-iwa}, we have
\[@
\rKp(t_{i-1},t_i)&=N_i \rKp^M\bar N_i,\\
\rKp(t_i)&=N_{i}\rKp^M\bar N_{i+1},\\
\rK(t_i)&=N_{i+1}\rK^M\bar N_i.
\]

Since \(\{\iota\}=\{\iota\}_{0\gamma}\) is \(\vec s\)-generic, the value
\(t=0\) lies on one of the open intervals \((t_{i},t_{i+1})\).  We may
and do asssume that \(0 \in (t_0,t_1)\).
Now put \(v_1=v\), and for \(i \geq 1\), define inductively
\[
v_{i+1}=\int_{N_{i+1}} x.v_{i}\,dx.
\]
We make two claims: (i) \(v_i\) is a non-zero multiple of
\(\int_{N_i}x.v\, dx\), and (ii) \(v_i\) is
\(\bigl(\rKp(t_{i-1},t_i),\theta(t_{i-1},t_i)\bigr)\)-isotypic, where \(\theta(t_{i-1},t_i)\) is the 
character of \(\rKp(t_{i-1},t_i)\) obtained by using the construction
of \ref{linear-setting} but with the embeddings
\(\{\iota\}_{t\gamma}\), \(t \in (t_{i-1},t_i)\) in place of \(\{\iota\}\).

We prove the claims by induction.  
(i) is a simple consequence of Fubini's theorem and the unimodularity
of \(U\).  

We have seen that \(\rKp(t_i,t_{i+1})=N_{i+1} \rKp^M \bar N_{i+1}\).
Since \(N_{i+1}, \bar N_{i+1} \subset
\ker\bigl(\theta(t_i,t_{i+1})\bigr)\) by Lemma~\ref{tri-on-nil}, it
follows that \(\ker\bigl(\theta(t_i,t_{i+1})\bigr)=N_{i+1} \ker(\theta_M) \bar N_{i+1}\).  In
particular, \(N_{i+1} \ker(\theta_M) \bar N_{i+1}\) is a subgroup.  Now we
prove (ii).  By the induction hypothesis, \(v_{i}\) is fixed by
\(\ker(\theta_M) \bar N_i\supset \ker(\theta_M) \bar N_{i+1}\).  Since
\(N_{i+1}\ker(\theta_M)\bar N_{i+1}\) is a compact subgroup, it
follows that \(v_{i+1}\) is fixed by this compact subgroup.  To finish
the proof of (ii), it suffices to show that \(\rKp^M\) acts on
\(v_{i+1}\) via the the character \(\theta_M\).  Indeed, for \(g \in
\rKp^M\),
\[
g.v_{i+1}=\int_{N_{i+1}} gxg^{-1}.g.v_i\,dx
=\int_{N_{i+1}} (gxg^{-1}).\theta_M(g)v_i\,dx=\theta_M(g)v_{i+1}.
\]
The last equality is because the compact group \(\rKp^M\) normalizes
\(N_{i+1}\).

It remains to prove the most important statement:

\begin{lemma*}
  \(v_{i+1} \neq 0\) for all \(i \geq 0\).
\end{lemma*}

\proof We prove this by induction on \(i\).  The case \(i=0\) holds by
assumption.  Let \((\tilde J,\tilde J_+,\tilde\varphi)\) be the triple
constructed in (\ref{h-new}) with \((M',M,G',G,\phi,r)=(M^{d-1}, M,
G^{d-1},G,\phi_{d-1},r_{d-1})\) and the embeddings
\(\{\iota\}_{t_i\gamma}\).  By Lemma~\ref{h-new} (d), \(N_{i+1} =
\rK(t_i)\cap U \subset (G',G)_{y,(0,s)} \cap U\) normalizes \(\tilde J
\cap U\), where \((G',G)_{y,(0,s)}\) is formed using the embeddings
\(\{\iota\}_{t_i\gamma}\).  It follows that \((N_{i+1} \cap G')(\tilde
J\cap U)\) is a subgroup, and in fact
\[
(N_{i+1} \cap G')(\tilde J\cap U) = N_{i+1}
\]
by [Yu, Lemmas 13.3 and 13.4].

It follows that \(v_{i+1}\) is a
non-zero multiple of
\[
\int_{N_{i+1} \cap G'} \int_{\tilde J \cap U}y.x.v_i \, dx\,dy
\]
We make two more claims: (iii) \(v'=\int_{\tilde J \cap U}x.v_i\,dx\)
is non-zero; (iv) \(v'\) is \(\bigl(\rKp(t_{i-1},t_{i})\cap
G',\theta(t_{i-1},t_i)\bigr)\)-isotypic.

To prove (iii), let \(W\) be the \((\tilde J_+,\tilde\varphi)\)-isotypic
subspace of \(V\).  This is naturally a representation of \(\tilde
J\).  It is easy to verify that \(\tilde J_+ \subset
\rKp(t_{i-1},t_i)\) and \(\theta(t_{i-1},t_i)|\tilde J_+=\tilde\varphi\).
Therefore, \(v_i \in W\) and 
the integral defining \(v'\) can be calculated within the \(\tilde
J\)-representation \(W\).
The image of \(\bar N_i \cap \tilde J\) in \(\tilde J/\ker(\tilde\varphi)\) is the same as
that of \(\tilde J^{\bar U}\), and the image of \(N_{i+1} \cap \tilde J\) in
\(\tilde J/\ker(\tilde\varphi)\) is the same as that of \(\tilde J^U\),
where \(\tilde J^U\) and \(\tilde J^{\bar U}\) are as in \ref{h-new}.
By (ii), \(\bar N_{i}\) fixes \(v_i\).  This implies
\(v' \neq 0\) by Lemmas~\ref{h-new} and \ref{L-nonzero}.

The proof of (iv) is similar to the proof of (ii) above.

Finally, we prove the lemma.  We can regard \(V\) as a smooth
representation of \(G'\).  We may also regard \(v'\) as a vector in
the representation \(V \otimes \phi_{d-1}^{-1}\) of \(G'\).  The
integral \(\int_{N_{i+1}\cap G'} y.v' \, dy\) is the same whether we
use the action of \(G'\) on \(V\) or on \(V \otimes \phi_{d-1}^{-1}\),
since \(\phi_{d-1}(N_{i+1} \cap G')=\{1\}\).  

We can now apply the setting of this section
to \(G'=G^{d-1}\), the twisted Levi
sequence \((G^0,\ldots,G^{d-1})\), the characters
\((\phi_0,\ldots,\phi_{d-2})\), and the
\((0,r_0/2,\ldots,r_{d-2}/2)\)-generic embeddings \(\{\iota\}_{t\gamma
}\), where \(t \in (t_i,t_{i+1})\).  Then we form the group \(\rKp'\)
and \(\theta '\) as in \ref{linear-setting}, and these are nothing but
\(\rKp'=\rKp(t_{i-1},t_i) \cap G'\), and \(\theta
'=(\theta(t_{i-1},t_i)|\rKp') \otimes \phi_{d-1}^{-1}\).
Recall that \(v' \in V \otimes \phi_{d-1}^{-1}\) is \((\rKp',\theta ')\)-isotypic by (iv).
Therefore, we can apply the induction hypothesis (for the induction on
\(d\)) to conclude \(v_{i+1} \neq 0\).  This proves the lemma,
and hence the theorem.\qed

\begin{cor} \label{C-linear}
  Let \(K\) be a compact open subgroup of \(G\) and \(\rho\)
  an irreducible smooth representation of \(K\).  Suppose that
  \((K,\rho)\) satisfies conditions {\rm (i)} and {\rm (ii)} in
  (\ref{D-cover}), \(K \supset \rKp\), and \(\rho|\rKp\) is
  \(\theta\)-isotypic.  Then \((K,\rho)\) is a \(G\)-cover of \((K\cap
  M, \rho|K\cap M)\).
\end{cor}

\proof 
For an irreducible smooth representation $V$ of $G$, let $V^\rho$ (resp. $V^\theta$) denote the $\rho$-isotypic (resp. $\theta$-isotypic) component in $V$.
Then, $V^\rho\subset V^\theta$. Since the natural map $V\rightarrow V_U$ on $V^\theta$ is injective, it is also injective on $V^\rho$. 
\qed

\section{Construction of types}\label{S-type}

{\it From now on we assume that \(p\) is odd.}

\subsection[Depth-zero datum] \label{depth0-datum}
We now review the construction of types
of depth zero by [MP2].   
We define a depth-zero datum to be a triple
\(\bigl((G,M),(y,\iota),(K_M,\rho_M)\bigr)\) such that
\begin{itemize}
\item[(i)] \(G\) is a connected reductive group over \(F\) and \(M\) a Levi
  subgroup of \(G\).
\item[(ii)] \(y \in \scrB(M)\) is such that \(M_{y,0}\) is a maximal
  parahoric subgroup of \(M\), and \(\iota:\scrB(M) \injto \scrB(G)\)
  is a \(0\)-generic embedding relative to \(y\).
\item[(iii)] \(K_M\) is a compact open subgroup of \(M\) containing
  \(M_{y,0}\) as a normal subgroup, and \(\rho_M\) is an irreducible smooth representation of
  \(K_M\) such that 
  \(\rho_M|M_{y,0}\) contains a cuspidal representation of \(M_{y,0:0+}\).
\end{itemize}

\begin{rem*} 
  Since \(M_{y,0}\) is a normal subgroup of \(K_M\), the restriction
  \(\rho_M|M_{y,0}\) is \(\rho '\)-isotypic where \(\rho '\) is any irreducible
  cuspidal representation occurring in \(\rho_M|M_{y,0}\).  It follows
  that \(\rho_M|M_{y,0+}\) is trivial (\({\bf 1}\)-isotypic).
\end{rem*}

This datum encodes not only a type of depth \(0\) in \(G\) but also how it arises
from a cover, as follows.
By [MP6, Proposition 6.8] and [BK, Proposition 5.4], \((K_{M},
\rho_{M})\) is an \(\grS\)-type where \(\grS\) is a finite set of the
form \(\{[M,\pi_1],\ldots,[M,\pi_n]\}\) with the \(\pi_i\)'s irreducible
supercuspidal representations of \(M\).  Note that {\it when \(K_{M}\) is the maximal
compact subgroup fixing \(y\) of $M$, \(\grS\) is a singleton.}

By Proposition~\ref{P-iwa} (b), \(K_{G}:=K_{M}
G_{\iota(y),0}\) is a subgroup such that 
\[
K_G/G_{\iota(y),0+} \simeq K_{M}/M_{y,0+}.
\]
Let \(\rho_G\) be the representation of \(K_G\)
obtained by composing the above isomorphism with \(\rho_{M}\).
Then \((K_{G},\rho_{G})\)  is a \(G\)-cover of \((
K_{M}, \rho_{M})\).  Therefore, by \cite[Theorem 8.3]{BK}, it is an \(\grS(G)\)-type, where \(\grS(G)\) is defined in [BK, \S8].

\subsection[The datum] \label{datum}
The datum $\Sigma$ from which we will construct a type is a
\(5\)-tuple 
\[
\Sigma:=\bigl((\vec G,M^0),(y,\{\iota\}),\vec
r,(K_{M^0},\rho_{M^0}),\vec\phi\bigr)
\]
entirely analogous to that in [Yu, \S3], as follows:
\begin{itemize}
\item [\bf D1] \(\vec G=(G^0,\ldots,G^d)\) is a tamely ramified
  twisted Levi sequence in \(G\), and \(M^0\) a Levi subgroup of
  \(G^0\).  Unlike [Yu, \S3], we impose {\bf no} assumption on
  \(\scrZ(G^0)/\scrZ(G)\). We
  construct a Levi subgroup \(M\) of \(G\) and a generalized twisted
  Levi sequence \(\vec M\) in \(M\) as in \ref{levi-seq}.
\item [\bf D2] \(y\) is a point in \(\scrB(M^0)\), and \(\{\iota\}\)
  is a commutative diagram of embeddings of buildings as in (\ref{emb}),
  \(\vec s\)-generic relative to \(y\), where \(\vec
  s=(0,r_0/2,\ldots,r_{d-1}/2)\).

\item [\bf D3] \(\vec r=(r_0,\ldots,r_d)\) is a sequence of real
  numbers satisfying \(0 < r_0 < r_1 <\cdots< r_{d-1}\leq r_d\) if \(d
  > 0\), \(0 \leq r_0\) if \(d=0\).
\item [\bf D4] 
  \(\bigl(G^0,M^0,y,\iota:\scrB(M^0)\injto \scrB(G^0),(K_{M^0},\rho_{M^0})\bigr)\)
  is a depth zero datum.
\item [\bf D5] \(\vec\phi=(\phi_0,\ldots,\phi_d)\), where \(\phi_i\) is a quasi-character of
  \(G^i\) such that \(\phi_i\) is \(G^{i+1}\)-generic of depth \(r_i\)
  relative to \(x\) for all \(x \in \scrB(G^i)\).
\end{itemize}

\begin{rem}  Again, the datum encodes not just the type itself but
  also how the type arises as a cover.
  In practice, one may start with a
  \(5\)-tuple \((\vec G,y,\vec
  r,\rho,\vec\phi)\) similar to [Yu, \S3], but with no assumption on
  \(\scrZ(G^0)/\scrZ(G)\), and instead of {\bf D4} of [Yu, \S3], we
  assume that \((G^0_{y,0},\rho)\) is an unrefined minimal \(K\)-type
  of depth \(0\) in the sense of Moy-Prasad.  We then
  construct \(M^0\) and \(\{\iota\}\) 
  as follows.  By [MP2, 6.3], to the parahoric subgroup \(G^0_{y,0}\) of
  \(G^0\), we can attach a Levi subgroup \(M^0\) of \(G^0\), unique up
  to conjugacy by \(G^0_{y,0}\).  From the construction there, we see
  that there is an embedding \(\iota:\scrB(M^0) \injto \scrB(G^0)\)
  whose image contains \(y\), \(M^0_{y,0}\) is maximal parahoric, and
  \(\iota\) is \(0\)-generic relative to \(y\).  One can then extend/modify
  \(\iota\) to a family \(\{\iota\}\) which is \(\vec s\)-generic by Lemma \ref{lem: genemb}.
  
  Of course, there are choices involved in this procedure.  Also,
  in order to get the finest $\grS$-types, i.e., \(\grS\)-types with
  \(\grS=\{\grs\}\) a singleton, we need to refine the datum
  \(\rho\) a little bit.  Eventually we end up with a datum as defined
  above.
\end{rem}

\subsection[The construction]\label{subsec: construction} We now put
\(K^0=K_{G^0}=K_{M^0}G^0_{y,0}=K_{M^0}\rK^0\) and \(\rho=\rho_{G^0}\)
as in (\ref{depth0-datum}), and for $i\geq 1$, put
\[
K^i=K^0G^1_{y,{s_0}}\cdots G^i_{y,s_{i-1}}=K_{M^0}\rK^i, \qquad K^i_+=\rKp^i.
\]
Again we remind the reader that these groups
depend on the choice of \(\{\iota\}\), and this \(K^i\) may be
different from the one used in [Yu, \S3 and \S4]. 
Nevertheless, it is easy to see that 
the construction in [Yu, \S4] can be carried out literally without any
modification to give a representation \(\rho_i\) for each \(K^i\), \(i
\geq 0\), with \(\rho_0=\rho\).

Moreover, \(\Sigma_M:=(\vec M,y,\vec r,\rho_M,\vec\phi)\) is a datum for
constructing a supercuspidal type in \(M\); see [Yu, Remark 15.4] and
the discussions following [Yu, Theorem 15.7].  So we can construct
supercuspidal types \((K^i_M, \rho_i^M)\), for each \(i\geq0\), where
\(K^0_M=K_{M^0}\),
\[
K_M^i=K_{M^0}M^1_{y,s_0}\cdots M^i_{y,s_{i-1}}, \qquad
i\geq 1.
\]
Write
\[
\scrT^i:=(K^i,\rho^i),\ \scrT:=\scrT^d=(K^d,\rho^d);\qquad \scrT^i_M:=(K^i_M,\rho_i^M),\ \scrT_M:=\scrT^d_M=(K^d_M,\rho_d^M).
\]
Let \(\grS_i\) be the finite set such that $\scrT^i_M$ is an
\(\grS_i\)-type in \(M^i\).  {\it If $K_{M^0}$ is the fixer of $y$ in $M^0$, $\grS_i$ is a singleton.}

For $\pi\in \scrR(G^i)$, we write $\scrT^i<\pi$ if $\rho^i$ occurs in $\pi|K^i$.

\begin{thm}\label{thm: types} For \(i \geq 0\), $\scrT^i$ is a \(G^i\)-cover of
$\scrT^i_M$. Hence $\scrT^i$ is an $\grS_i(G^i)$-type in \(G^i\).
\end{thm}

\proof The second statement follows from the first and \cite[Theorem 8.3]{BK}.

Condition (i) in \ref{D-cover} follows from
Proposition~\ref{P-iwa} (b).

We now verify condition (ii) in \ref{D-cover} by induction.  The case
of \(i=0\) is just the definition.  The inductive construction of
\(\rho_i\) from \(\rho_{i-1}\) in [Yu, \S4] relies on the Heisenberg
triple \((J^i,J^i_+,\varphi_i)\), where
\(\varphi_i=\hat\phi_{i-1}|J^i_+\).  Similarly, to construct
\(\rho_i^M\) from \(\rho_{i-1}^M\) we use an analogous Heisenberg
triple \((J^i_M, J^i_{M,+},\varphi_i^M)\).  It follows from the
definitions of these objects that \(J^i_M = J^i \cap M_i\),
\(J^i_{M,+} = J^i_+ \cap M_i\), and \(\varphi_i^M =
\varphi_i|J^i_{M,+}\).  Moreover, \(J^i\) and \(J^i_+\) decompose
with respect to \((U^i,M^i,\bar U^i)\) by Proposition~\ref{P-iwa},
where \(U^i=U \cap G^i\), \(\bar U^i=\bar U \cap G^i\).  Since
\(\{\iota\}\) is \(s\)-generic, we have \(J^i \cap U^i = J^i_+ \cap
U^i\) and \(J^i \cap \bar U^i=J^i_+ \cap \bar U^i\).  It follows that
the inclusion \(J^i_M \subset J^i\) induces an isomorphism
\[
J^i_M/J^i_{M,+}\simeq J^i/J^i_+.
\]
Let \(N^i=\ker(\varphi_i)\), \(N^i_M=\ker(\varphi_i^M)\).  We can verify
that the following diagram is commutative: 
\[
\XYmatrix{
J_M^i/N_M^i \xyto[d]^{\wr} \xyto^{\sim} &(J_M^i/J^i_{M,+})^\sharp\xyto[d]^{\wr}\\
J^i/N^i \xyto^{\sim} &(J^i/J^i_+)^\sharp,\\
}
\]
where  $(J_M^i/J^i_{M,+})^\sharp$ and $(J^i/J^i_{+})^\sharp$ are defined as in \cite[\S10]{Yu}.
The vertical arrows are the isomorphisms induced by inclusion,
and the horizontal arrows are the canonical special isomorphisms
constructed in [Yu, Proposition 11.4].  In addition, the following
diagram is also commutative:
\[
\XYmatrix{
K_M^{i-1} \xysubset[d] \xyto & \Sp(J_M^i/N_M^i)\xyto[d]\\
K^{i-1} \xyto &\Sp(J^i/N^i),
}
\]
where the horizonal arrows are induced by conjugations, and the
vertical arrow on the right is induced by the isomorphism
\(J_M^i/J_{M,+}^i \simeq J^i/J^i_+\).  It follows from these and the
definitions of \(\rho_i\) and \(\rho_i^M\) that we do have \(\rho_i|K^i_M=\rho_i^M\).

Define \((K^i_+,\theta_i)\) as in the preceding section.  Note that this
\(K^i_+\) is identical to the \(K^i_+\) used in [Yu].  By [Yu,
Proposition 4.4], \(\rho_i|K^i_+\) is \(\theta_i\)-isotypic.
By Proposition~\ref{tri-on-nil}, \(K^i_+ \cap U^i, K^i_+ \cap \bar U^i
\subset \ker(\rho_i)\).  Since \(K^i_+ \cap U^i = K^i \cap U^i\) and
\(K^i_+ \cap \bar U^i=K^i \cap \bar U^i\) by the genericity of $\{\iota\}$, we have proved condition (ii)
of \ref{D-cover} completely.

Now we see that all hypotheses for Corollary~\ref{C-linear} are satisfied.
The theorem is proved.\qed

\section{Support of Hecke algebras}

Let \((K^i,\rho_i)\) be as in \S\ref{subsec: construction}. Let \(\check\rho_i\) be the contragradient of \(\rho_i\).
Then the Hecke algebra \(\Hecke(G^i,\rho_i)\) associated to \((K^i,\rho_i)\) is
defined as follows:
\[
\Hecke(G^i,\rho_i)=\{f\in C_c(G^i,\End(\check\rho_i))\mid f(jgj')=\check\rho_i(j)f(g)\check\rho_i(j')\ \textrm{for }g\in G^i,\ j,j'\in K^i\}.
\]
As in \cite[\S17]{Yu}, we write \(\check\Hecke(G^i,\rho_i)\) for \(\Hecke(G^i,\check\rho_i)\).
For $g\in G^i$, let \(I_g(\rho_i)\) denote the space of intertwining maps
\(\textrm{Hom}_{K^i\cap gK^ig^{-1}}(\,^g\!\rho_i,\rho_i)\) where $^g\!\rho_i$ is a representation of $gK^ig^{-1}$ with $^g\!\rho_i(h)=\rho_i(g^{-1}hg)$ for $h\in gK^ig^{-1}$ (see also \cite[p582]{Yu}).

\begin{thm}\ 
\begin{itemize}
\item[(a)]
The support of \(\check\Hecke(G^i,\rho_i)\) is contained in \(K^iG^0K^i\).
\item[(b)]
For $g\in G^0$, we have
\[
I_g(\rho_i)=I_g(\rho_0\mid K^0)\otimes I_g(\tilde\phi_0)\otimes\cdots\otimes I_g(\tilde\phi_{i-1})
\]
where $I_g(\tilde\phi_j)$ is 1-dimensional for $j=0,\cdots,d-1$.
\end{itemize}
\end{thm}

Again the proof  in \cite[\S15]{Yu} can be carried out without change.

\begin{cor}
The support of \(\check\Hecke(G^i,\rho_i)\) is contained in \(K^iN_{G^0}(M^0)K^i\).
\end{cor}

\proof
By \cite[Theorem 4.15]{Mo}, for $g\in G^0$, we have $I_g(\rho_0\mid K^0)\ne 0$ only if $g\in K^0 N_{G^0}(M^0) K^0$.
Hence, combining with the above theorem, the corollary follows.
\qed

\section{Exhaustion}\label{sec: exhaustion}

Recall that in \cite{Kim}, it is proved that all supercuspidal representations arise from the construction given in \cite{Yu} under some hypotheses (see \cite[\S3.4]{Kim}). In this section, we prove that our construction gives all types parameterizing $\scrR^\grs$, $\grs\in\scrI$ (see (1) in \S1) under the same hypotheses $\mathrm{(H{\mathit k}), (HB), (HGT)}$ and $\mathrm{(H{\scrN})}$ (see \cite[\S3.4]{Kim} for details). We adopt notation and terminologies from \cite{Kim}.

\begin{thm} Suppose $\mathrm{(H{\mathit k}), (HB), (HGT)}$ and $\mathrm{(H{\scrN})}$ are valid. 
For each inertial class $\grs\in\scrI$, there is a datum $\bigl((\vec G,M^0),(y,\{\iota\}),\vec
r,(K_{M^0},\rho_{M^0}),\vec\phi\bigr)$ so that $(K^d, \rho_d)$ is an $\grs$-type.
\end{thm}

\noindent{\sc Sketch of the proof.} 
Let $\scrE^t(G)$ be the set of irreducible smooth tempered representations. Note that for any $\grs=[(M_{\grs},\pi_{\grs})] \in\scrI$, the Plancherel measure of $\scrE^t(G)\cap\scrR^\grs$ is nonzero.
Hence, it is enough to show that there are $\pi \in \scrE^t(G)\cap\scrR^\grs$ and a datum $\bigl((\vec G,M^0),(y,\{\iota\}),\vec r,(K_{M^0},\rho_{M^0}),\vec\phi\bigr)$ so that $(K^d,\rho_d)$ gives a $G$-cover of the supercuspidal type of the cuspidal pair $(M_{\grs},\pi_{\grs})$ which supports $\pi$. We sketch a proof in several steps below:

\vspace{3pt}

(1) Recall from \cite[\S4]{KM}, for a given strongly good datum $(\vec G,x,\vec r,\vec\phi)$, upon fixing embeddings of buildings, $\scrB(G^0)\hookrightarrow\scrB(G^1)\hookrightarrow\cdots\hookrightarrow\scrB(G^d)$, one can associate a strongly good type $(K^d_{x+},\theta_d)$ where $K^d_{x+}=G^0_{x,0+}G^1_{x,s_0+}\cdots G^d_{x,s_{d-1}+}$ and $\theta_d$ is constructed as in \ref{linear-setting} (see \cite[\S4]{KM} for details). 
Then, by \cite[Theorem 11.4]{Kim}, there are $\pi\in\scrE^t(G)\cap\scrR^\grs$ and a strongly good type $(K^d_{x+},\theta_d)$ such that $(K^d_{x+},\theta_d)<\pi$.

(2) Let $V_\pi^{\theta_d}$ be the $\theta_d$ isotypic component of $V_\pi$. Then, $V_\pi^{\theta_d}$ is stabilized by $G^0_{x,0}$. Let $y\in\scrB(G^0)$ be such that (i) $G^0_{y,0}\subset G^0_{x,0}$, (ii) $V_\pi^{\theta_d}$ has nontrivial $G^0_{y,0^+}$-invariants, and (iii) $G^0_{y,0}$ is minimal satisfying (i) and (ii). Such a $y$ exist since $\theta_d$ is trivial on $G^0_{x,0^+}$.


(3) Let $M^0$ be a Levi subgroup of $G^0$ so that $M^0_{y,0}$ is a maximal parahoric subgroup of $M^0$, and $M$ the Levi subgroup given by the centralizer of $\scrZ_\rms(M^0)$. Let $P=MU$ be the parabolic subgroup of $G$ so that $(P\cap G^0_{x,0})G^0_{y,0+}=G^0_{y,0}$. Let $\overline U$ be the opposite unipotent radical.

(4) Form $K_\vdash=(\overline UM\cap K^d_{x+})(U\cap (G^0_{y,0+}G^1_{x,s_0}\cdots G^d_{x,s_{d-1}}))$, which is defined in \cite[\S13]{Kim}. Note that $K_{x+}\subset K_\vdash$ and $\vec\phi$ defines a character $\theta_d'$ of $K_\vdash$ such that $\theta'_d|K^d_{x+}=\theta_d$. 
Then, by \cite[Corollary 13.12]{Kim}, $(K_\vdash,\theta_d')<V_\pi$, and $V_\pi^{\theta_d'}\subset V_\pi^{\theta_d}$. 

(5) Consider $V_\pi^{\theta_d'}$. Note that $V_\pi^{\theta_d'}$ is stabilized by $K_{y,M}=M^0_{[y]} (M\cap (G^1_{y,s_0}\cdots G^d_{y,s_{d-1}}))=M^0_{[y]}M^1_{y,s_0}\cdots M^d_{y,s_{d-1}}$ where $M^0_{[y]}$ is the stabilizer of the image $[y]$ of $y$ in the reduced building of $M$ and $M^i=M\cap G^i$. Let $\hat\theta_d'$ be an irreducible representation of  $K_{y,M}$ such that $\hat\theta'_d|K_\vdash$ is $\theta'_d$-isotypic as in \cite[\S13]{Kim} ($\hat\theta_d'$ is denoted by $\kappa$ in \cite[\S13]{Kim}). Then, by \cite[Corollary 18.6]{Kim}, there is an irreducible representation $\tau'$ of $M^0_{[y]}$ factoring through $Z_MM^0_{y,0^+}$ such that $\tau'\otimes\hat\theta_d'$ is contained in $V_\pi^{\theta_d'}$. Then, $\tau'|M^0_{y,0}$ induces a cuspidal representation of $M^0_{y,0}/M^0_{y,0^+}$ since otherwise, there is a smaller parahoric subgroup $G^0_{z,0}\subset G^0_{x,0}$ with nontrivial $G^0_{z,0+}$ invariants in $V_\pi^{\theta_d}$, which is a contradiction to the choice of $G^0_{y,0}$ in (2). 

(6) From the proof of \cite[Theorem 19.1]{Kim}, $\pi_M:=c\textrm{-}\mathrm{Ind}_{K_{y,M}}^M\tau'\otimes\hat\theta'_d$ is a supercuspidal representation of $M$ associated to a generic datum $(\vec M,y,\vec r,\vec\phi_M,\tau')$ where $\vec M=(M^0,M^1,\cdots, M^d)$ and $\vec\phi_M=(\phi_0|M,\cdots,\phi_d|M)$. Moreover, 
$(M, \pi_M)$ is equivalent to $\grs$ in $\scrI$.

(7) Let $K_{M^0}$ be the maximal compact subgroup of $M^0_{[y]}$ and $\rho_{M^0}$ a subrepresentation of $\tau'$ when restricted to $K_{M^0}$. Then, $((\vec M,M^0),y,\vec r,(K_{M^0},\rho_{M^0}),\vec\phi_M)$ gives a supercuspidal type $(K^d_M,\rho^M_d)$.

(8) Consider $((\vec G,M^0), (y,\{\iota\}),\vec r,(K_{M^0},\rho_{M^0}),\vec\phi)$ where $\{\iota\}$ is $(0,s_0,\cdots,s_{d-1})$-generic. Then, by construction, $(K^d,\rho_d)$ is a cover of $(K_M^d,\rho^d_M)$, hence an $\grs$-type.
\qed

\begin{rem*}
\begin{enumerate}
\item The above proof starts with $\pi\in\scrR^{\grs}$ and a strongly good type contained in $\pi$, proceeds  toward nailing down a supercuspidal type $(K_M, \rho_{d}^M)$ out of the strongly good type, and then finally finds a type as a cover $(K^d,\rho_d)$ of $(K_M,\rho_{d}^M)$. On the other hand, we note that it is possible to start with a supercuspidal type datum for $\grs$ and work toward getting a cover. However, to nail down $\vec G$ in the datum, we find the proof above more efficient. 

\item A priori, we can not assume our choices of $y$ in (2) or embeddings $\scrB(G^0)\hookrightarrow\scrB(G^1)\hookrightarrow\cdots\hookrightarrow\scrB(G^d)$ in (1) satisfy any genericity condition. Hence, we still need to work with an auxiliary group $K_\vdash$ in (4) (cf. Remark \ref{rmk: vdash}).
\end{enumerate}
\end{rem*}

\section{Equivalence}\label{sec: equivalence}

\begin{defn}\label{defn: equiv} Let $\Sigma$ and $\dot\Sigma$ be two data as in \S\ref{datum}. Let $\scrT=(K,\rho)$ (resp. $\dot\scrT$) be the type constructed in \S\ref{subsec: construction}  associated to $\Sigma$ (resp. $\dot\Sigma$).

\begin{enumerate}
\item[(i)]
Define $\scrR_{\scrT}$ to be the category of smooth representations $\pi$ which are generated by the $\rho$ isotypic components of $V_\pi$.

\item[(ii)] Let $\dot\scrT$ be the type associated to $\dot\Sigma$. We say that $\scrT$ and $\dot\scrT$ are equivalent if there is $\grS\subset\scrI$ such that $\scrR^{\grS}=\scrR_{\scrT}=\scrR_{\dot\scrT}$ where $\scrR^{\grS}=\prod_{\grs\in\grS}\scrR^{\grs}$.

\end{enumerate}
\end{defn}

From now on, we assume that our data $\Sigma$ satisfy the hypothesis $\mathrm{C}(\vec{G})$ in \cite[\S2.6]{HaMu}.

\begin{thm}\label{thm: sc equiv} Let $\Sigma:=\bigl((\vec G,G^{0}),(y,\{\iota\}),\vec
r,(K_{G^{0}},\rho_{G^{0}}),\vec\phi\bigr)$ and 
$\dot\Sigma:=\bigl((\vec{\dot G},\dot G^{0}),(\dot y,\{\dot\iota\}),\vec\dot r,(K_{\dot G^{0}},\rho_{\dot G^{0}}),\vec\phi\bigr)$ be two data such that $Z_{G^{0}}/Z_G$ (resp. $Z_{\dot G^0}/Z_G$) is $F$-anisotropic and $K_{G^{0}}$ (resp. $K_{\dot G^{0}}$) is the maximal compact subgroup of $G^{0}_{[y]}$ (resp. $\dot G^{0}_{[\dot y]}$). Let $\scrT:=(K,\rho)$ and $\dot\scrT:=(\dot K,\dot\rho)$.
Let $\phi=\prod_{i=0}^{d}\phi_i$ and $\dot\phi=\prod_{i=0}^{\dot d}\dot\phi_i$ be characters of $G^0$ and $\dot G^0$ respectively. Then, $\scrT$ and $\dot\scrT$ are equivalent if and only if there is $g\in G$ such that $gy=\dot y$, $^g\!K_{G^{0}}=K_{\dot G^{0}}$ and $^g\!(\rho_{G^{0}}\otimes\phi)\simeq (\dot\rho_{\dot G^{0}}\otimes\dot\phi)$ as $K_{\dot G^0}$ representations.
\end{thm}

\proof
Note that $\scrT$ and $\dot\scrT$ are supercuspidal types. 

Suppose $\scrR_{\scrT}=\scrR_{\dot\scrT}$. Let $\pi\in\scrR_{\scrT}=\scrR_{\dot\scrT}$ be irreducible supercuspidal. Then, let $\tilde\rho_{G^0}$ (resp. $\tilde{\dot\rho}_{\dot G^0}$) be a representation of $G^0_{[y]}$ (resp. $\dot G^0_{[\dot y]}$) containing $\rho_{G^0}$ (resp. $\dot\rho_{\dot G^0}$) so that $(\vec G, y, \vec r,\vec\phi,\tilde\rho_{G^0})$ (resp. $(\vec {\dot G},\dot y, \vec{\dot r},\vec{\dot\phi},\tilde{\dot\rho}_{\dot G^0})$) is a supercuspidal datum for $\pi$. Then, \cite[Theorem 6.7]{HaMu}, there is $g\in G$ so that $gy=\dot y$, $^g\!G^0=\dot G^0$, and $^g(\tilde\rho_{G^0}\otimes\phi)\simeq(\tilde\rho_{\dot G^0}\otimes\dot\phi)$ as $\dot G^0_{[\dot y]}$ representations. Since both $^g\!(\rho_{G^0}\otimes\phi)$ and  $(\dot\rho_{\dot G^0}\otimes\dot\phi)$ are subrepresentations of $(\tilde\rho_{\dot G^0}\otimes\dot\phi)$ and $K_{\dot G^0}$ is a normal subgroup of $\dot G^0_{[y]}$, there is a $\dot g\in \dot G^0_{[\dot y]}$ so that $^{\dot gg}\!(\rho_{G^{0}}\otimes\phi)\simeq (\dot\rho_{\dot G^{0}}\otimes\dot\phi)$. 

Conversely, suppose there is $g\in G$ such that $gy=\dot y$, $^g\!K_{G^{0}}=K_{\dot G^{0}}$ and $^g\!(\rho_{G^{0}}\otimes\phi)\simeq (\dot\rho_{\dot G^{0}}\otimes\dot\phi)$ as $K_{\dot G^0}$ representations. It is enough to show that there is a supercuspidal representation $\pi\in\scrR_{\scrT}\cap\scrR_{\dot\scrT}$. Let $\tilde\rho_{G^0}$ (resp. $\tilde{\dot\rho}_{\dot G^0}$) be a representation of $G^0_{[y]}$ (resp. $\dot G^0_{[\dot y]}$) containing $\rho_{G^0}$ (resp. $\dot\rho_{\dot G^0}$) so that $^g(\tilde\rho_{G^0}\otimes\phi)\simeq(\tilde\rho_{\dot G^0}\otimes\dot\phi)$ as $\dot G^0_{[\dot y]}$ representations. Then, $(\vec G, y, \vec r,\vec\phi,\tilde\rho_{G^0})$ and $(\vec {\dot G},\dot y, \vec{\dot r},\vec{\dot\phi},\tilde{\dot\rho}_{\dot G^0})$ are supercuspidal data, and their associated supercuspidal representations are isomorphic, which are in $\scrR_{\scrT}\cap\scrR_{\dot\scrT}$. Hence, $\scrR_{\scrT}=\scrR_{\dot\scrT}$.
\qed

\begin{thm} Let $\Sigma:=\bigl((\vec G,M^{0}),(y,\{\iota\}),\vec
r,(K_{M^{0}},\rho_{G^{0}}),\vec\phi\bigr)$ and 
$\dot\Sigma:=\bigl((\vec{\dot G},\dot M^{0}),(\dot y,\{\dot\iota\}),\vec{\dot r},(K_{\dot M^{0}},\rho_{\dot M^{0}}),\vec\phi\bigr)$ be two data such that $K_{M^{0}}$ (resp. $K_{\dot M^{0}}$) is the maximal compact subgroup of $M^{0}_{[y]}$ (resp. $\dot M^{0}_{[\dot y]}$). Let $\scrT:=(K,\rho)$ and $\dot\scrT:=(\dot K,\dot\rho)$.
Let $\phi=\prod_{i=0}^{d}(\phi_i|M^0)$ and $\dot\phi=\prod_{i=0}^{\dot d}(\dot\phi_i|\dot M^0)$ be characters of $M^0$ and $\dot M^0$ respectively. Then, $\scrT$ and $\dot\scrT$ are equivalent if and only if there is $g\in G$ such that $^g\!K_{M^0}=K_{\dot M^{0}}$ and $^g\!(\rho_{M^{0}}\otimes\phi)\simeq (\dot\rho_{\dot M^{0}}\otimes\dot\phi)$ as $K_{\dot M^0}$ representations.
\end{thm}

\proof 
Let $M$ (resp. $\dot M$) be the centralizer of $\scrZ_s(M^0)$ (resp. $\scrZ_s(\dot M^0)$) in $G$ as in \S\ref{levi-seq}. Let $\grs$ (resp. $\dot\grs$) be the inertial class of $(M,\pi_{M})$ (resp. $(\dot M,\dot\pi_{\dot M})$) where $\pi_{M}$ (resp. $\dot\pi_{\dot M}$) is the supercuspidal representation such that $\scrT_M<\pi_M$ (resp. $\dot\scrT_{\dot M}<\dot\pi_{\dot M}$). Then, we have $\scrR_{\scrT}=\scrR^{\grs}$ and $\scrR_{\dot\scrT}=\scrR^{\dot\grs}$.

Suppose $\scrR_{\scrT}=\scrR_{\dot\scrT}$, thus $\scrR^{\grs}=\scrR^{\dot\grs}$. Then, there is an unramified character $\chi$ of $M^0$ so that $^g\!M=\dot M$ and $^g\!\pi_M\simeq\dot\pi_{\dot M}\otimes\chi$. By Theorem \ref{thm: sc equiv}, there is $\dot m\in \dot M$ so that $^{\dot mg}\!(\rho_{M^0}\otimes\phi)\simeq \rho_{\dot M^0}\otimes(\dot\phi\chi)=\rho_{\dot M^0}\otimes\dot\phi$. Since $\chi$ is trivial on $\dot K$, we have $^{\dot mg}\!(\rho_{M^0}\otimes\phi)\simeq=\rho_{\dot M^0}\otimes\dot\phi$ as representations of $\dot K$.

Conversely, suppose there is $g\in G$ such that $^g\!K_{M^0}=K_{\dot M^{0}}$ and $^g\!(\rho_{M^{0}}\otimes\phi)\simeq (\dot\rho_{\dot M^{0}}\otimes\dot\phi)$. By Theorem \ref{thm: sc equiv}, $\scrR_{^g\!\scrT_M}=\scrR_{\dot\scrT_{\dot M}}\subset \scrR(M)$. Since $^g\!\scrT$ and $\dot\scrT$ are covers of $^g\!\scrT_M$ and $\dot\scrT_{\dot M}$ respectively, we have $\scrR_{\scrT}=\scrR_{^g\!\scrT}=\scrR_{\dot\scrT}$.
\qed

\begin{rem*} In \cite{Ka}, Kaletha studied the equivalence of regular representations. His methodology, especially in \S3.5, may allow the replacement of the hypothesis $C(\vec G)$ by a weaker one. 
\end{rem*}


\begin{thebibliography}{BDKV}
\bibitem[Be]{Be} J. Bernstein: {\it Le ``centre'' de Bernstein}, Edited by P. Deligne. Travaux en Cours, Representations of reductive groups over a local field, 1–32, Hermann, Paris, 1984.
\bibitem[Bl]{Bl} C. Blondel: {\it Crit\`ere d'injectivit\'e pour l'application de Jacquet}, C. R. Acad. Sci. Paris Sér. I Math. {\bf 325} (1997), no. 11, 1149 --1152.
\bibitem[Bo]{Bo} A. Borel: {\it Linear algebraic groups}, 2nd edition. Graduate Texts in Mathematics, {\bf126} Springer-Verlag, New York, 1991
\bibitem[BT1]{BT1} F. Bruhat and J. Tits: {\it Groupes r\'eductifs sur un corps local}, Inst. Hautes Études Sci. Publ. Math. {\bf 41} (1972), 5--251
\bibitem[BK]{BK} C.J. Bushnell and P.C. Kutzko: {\it Smooth representations of reductive
\(p\)-adic groups: structure theory via types}, Proc.\ London Math.\ Soc.~{\bf 77}  (1998),
582--634.
\bibitem[GR]{GR} D. Goldberg and A. Roche: {\it Types in $SL_n$}, Proc. London Math. Soc. (3) {\bf 85} (2002), no. 1, 119--138.
\bibitem[HaMu]{HaMu} J. Hakim and F. Murnaghan: {\it Distinguished tame supercuspidal representations.} Int. Math. Res. Pap. IMRP (2008),  no. 2

\bibitem[HM]{HM} R.Howe:  {\it Harish-Chandra homomorphisms for p-adic groups}, With the collaboration of A. Moy, CBMS Regional Conference Series in Mathematics,  {\bf 59}, Published for the Conference Board of the Mathematical Sciences, Washington, DC; by the American Mathematical Society, Providence, RI, 1985

\bibitem[Ka]{Ka} T. Kaletha: {\it Regular supercuspidal representations}, preprint, {\tt arXiv:1602.03144}

\bibitem[K]{Kim} J.-L. Kim: {\it Supercuspidal representations: an exhaustion theorem}, J. Amer. Math. Soc. {\bf 20} (2007), no. 2, 273--320

\bibitem[KM]{KM} J.-L. Kim and F. Murnaghan: {\it $K$-types and $\Gamma$-asymptotic expansions},  J. Reine Angew. Math. {\bf 592} (2006), 189--236

\bibitem[Mo]{Mo} L. Morris: {\it Tamely ramified intertwining algebras}, Invent. Math. {\bf 114} (1993), no. 1, 1–54.
\bibitem[MP1]{MP1} A. Moy and G. Prasad: {\it Unrefined minimal K-types for p-adic groups}, Invent. Math. {\bf 116} (1994), no. 1-3, 393--408.
\bibitem[MP2]{MP2} A. Moy and G. Prasad: {\it Jacquet functors and unrefined minimal K-types}, Comment. Math. Helv. {\bf 71} (1996), no. 1, 98--121.
\bibitem[SS]{SS} P. Schneider and U. Stuhler: {\it Representation theory and sheaves on the Bruhat-Tits building}, Inst. Hautes \'Etudes Sci. Publ. Math.  {\bf 85},   (1997) 97--191.
\bibitem[St]{St} R. Steinberg: {\it Torsion in reductive groups}, Adv.~in Math.~{\bf 15},
63--92 (1975).
\bibitem[Yu]{Yu} J.-K. Yu: {\it Construction of tame supercuspidal representations}, Journal of the A.M.S.~{\bf 14},  (2001) 579--622.
\end{thebibliography}
\end{document}